\documentclass[11pt]{amsart}
\usepackage{eucal}
\usepackage[all]{xy}
\usepackage{amsfonts,amssymb}
\usepackage{epsfig}
\usepackage[usenames]{color}
\usepackage{float} 
\xyoption{arc}

\newtheorem{theorem}{Theorem}[section]
\newtheorem{lemma}[theorem]{Lemma}

\newtheorem{corollary}[theorem]{Corollary}
\newtheorem{defn}[theorem]{Definition}

\newcommand{\Q}{\mathbb Q}
\newcommand{\C}{\mathbb C}
\newcommand{\R}{\mathbb R}

\newcommand{\D}{\mathbb D}

\newcommand{\Z}{\mathbb Z}
\newcommand{\N}{\mathbb N}
\newcommand{\J}{\mathcal J}
\newcommand{\M}{\mathcal M}
\newcommand{\K}{\mathbb K}

\title{A spectral sequence for symplectic homology}
\author{Mark McLean}

\begin{document}

\begin{abstract}
We construct a spectral sequence converging to symplectic homology
of a Lefschetz fibration whose $E^1$ page is related to Floer homology
of the monodromy symplectomorphism and its iterates.
We use this to show the existence of fixed points of certain symplectomorphisms.
\end{abstract}

\maketitle

\bibliographystyle{halpha}


\tableofcontents

\section{Introduction}

Symplectic homology is a useful tool in symplectic geometry.
It has been used in many areas.
We are interested in symplectic homology of Liouville domains.
A Liouville domain is
a compact manifold $F$ with boundary and a $1$-form $\lambda$
satisfying:
\begin{enumerate}
\item $d\lambda$ is a symplectic form.
\item $\alpha := \lambda|_{\partial F}$ is a contact form on the boundary of $F$.
\item The boundary is positively oriented. This means that if we have an outward
pointing vector field $V$ defined near $\partial F$, then the volume form:
\[i(V) (d\lambda)^n|_{\partial F} \]
has the same orientation as $\alpha \wedge (d\alpha)^n$.
\end{enumerate}
There exists a collar neighborhood $(1-\epsilon,1] \times \partial F$
of $\partial F$ such that $\lambda = r\alpha$ where
$r$ parameterizes $(1-\epsilon,1]$.
The reason why this exists is because we have a natural vector field
$K$ transverse to the boundary which is the $d\lambda$-dual of $\lambda$,
and the collar neighborhood is constructed by flowing the boundary $\partial F$
backwards along $K$.
We can form the completion of a Liouville domain by
extending $(1-\epsilon,1] \times \partial F$
by attaching a cylindrical end $[1,\infty) \times \partial F$
and extending $\lambda$ by $r\alpha$.
We write $\widehat{F}$ for the completion of a Liouville domain $F$.


Symplectic homology was a tool defined in \cite{Viterbo:functorsandcomputations}.
Viterbo used it to give obstructions to various Lagrangian embeddings in cotangent bundles,
and also used it to prove the Weinstein conjecture for subcritical Stein manifolds.
It also has many other applications.
Symplectic homology involves taking some time dependent Hamiltonian
$H : S^1 \times \widehat{F} \rightarrow \R$
such that $H = r^2$ near infinity, and creating a chain complex
involving fixed points of the time $1$ Hamiltonian symplectomorphism $\phi^1_H$.
We write $SH_*(\widehat{F})$ for this symplectic homology group.
See section \ref{section:symplectichomologydefn} for a precise definition.

Let $\phi : F \rightarrow F$ be a symplectomorphism which is the identity on the boundary $\partial F$.
We can assign a group $HF^*(\phi)$ to it called Floer cohomology.
This is basically a cohomology group whose chain complex is generated
by fixed points of $\phi$ away from the boundary.
In particular if all the fixed points of $\phi$ are non-degenerate, then we get
that the number of fixed points of $\phi$ is bounded below
by the rank of $HF^*(\phi)$.
The aim of this paper is to relate this group with symplectic homology.
We have to do this via another Floer homology group $HF_*(\phi,k)$.
The chain complex is generated by two copies of fixed points of $\phi^k$ modulo
an equivalence relation given by identifying a fixed point $x$
with $\phi^l(x)$ for some $l \in \Z$.

Suppose that $\pi : E \rightarrow \C$ is a Lefschetz fibration with one positive end whose
smooth fibers are symplectomorphic to $\widehat{F}$.
This is defined in section \ref{section:lefschetzfibrations}.
This has a compactly supported monodromy map $\phi' : \widehat{F} \rightarrow \widehat{F}$.
Suppose that $\phi'|_{F} = \phi$ and $\phi'$ is the identity outside $F \subset \widehat{F}$.
This Lefschetz fibration is symplectomorphic to $\widehat{M}$ for some Liouville domain $M$.

The first main theorem says:
\begin{theorem} \label{thm:spectralsequence}
If the dimension $2n$ of $E$ is greater than $2$ then
there is a spectral sequence converging to
$SH_*(\widehat{M})$ with $E^1$ page satisfying:
\[E^1_{0,q} = H^{n-q}(M)\]
and for $p > 1$,
\[E^1_{p,q} = HF_{q-p}(\phi,p).\]
For $p<0$,
\[E^1_{p,q} = 0.\]
The differential $d^r$ on the $r$th page sends
$E^r_{p,q}$ to $E^r_{p-r,q+r-1}$.
\end{theorem}
This is basically a slightly more refined statement
of \cite[Theorem 2.24]{McLean:symhomlef}.
Here is an example:
Suppose we have a trivial Lefschetz fibration
\[\pi : \widehat{F} \times \C \twoheadrightarrow \C\]
where $F$ is a Stein domain.
We have that the monodromy map is the identity map
and so 
\[HF_*(\phi,p) = H^{n-1-*}(\widehat{F}) \oplus
H^{n-*}(\widehat{F})\]
which gives us the following $E^1$ page:

\[
\begin{array}{c|cccc}
 & p=0 & p=1 & p=2 & \\
\hline
q = 3 & H^{n-3}(\widehat{F}) &
H^{n-3}(\widehat{F}) \oplus H^{n-2}(\widehat{F}) & 
H^{n-2}(\widehat{F}) \oplus H^{n-1}(\widehat{F}) & \cdots \\
q = 2 & H^{n-2}(\widehat{F}) &
H^{n-2}(\widehat{F}) \oplus H^{n-1}(\widehat{F}) &
H^{n-1}(\widehat{F}) \oplus H^n(\widehat{F}) & \\
q = 1 & H^{n-1}(\widehat{F}) &
H^{n-1}(\widehat{F}) \oplus H^n(\widehat{F}) & H^n(\widehat{F}) & \\
q = 0 & H^n(\widehat{F}) & H^n(\widehat{F}) & 0 & \\
q = -1 & 0 & 0 & 0 &  \\
q = -2 & 0 & 0 & 0 & \cdots \\
\vdots & & & \vdots &
\end{array}
\]
The differential $d^1_{p,q}$ has degree $(-1,0)$.
It is a natural projection map composed with a natural
inclusion map.
For instance it projects
\[E^1_{2,3} = H^{n-2}(\widehat{F}) \oplus H^{n-1}(\widehat{F})\]
to $H^{n-2}(\widehat{F})$ and then includes this into
\[E^1_{1,3} = H^{n-3}(\widehat{F}) \oplus H^{n-2}(\widehat{F})\]
via the natural inclusion.
It turns out that the spectral sequence in this case degenerates on the
first page giving us $SH_*(\widehat{F} \times \C) = 0$
basically by \cite[Theorem C]{Oancea:kunneth}.

We also have a relationship between $HF_*(\phi,k)$
and $HF_*(\phi^k,1)$:
\begin{theorem} \label{thm:groupactions}
Suppose that the coefficient field $\K$
has characteristic $0$ or characteristic $p$
where $p$ does not divide $k$,
then there exists a $\Z / k\Z$ action
$\Gamma$ on $HF_*(\phi^k,1)$
such that
$HF_*(\phi,k) \cong HF_*(\phi^k,1)^\Gamma$.
\end{theorem}
This is proven in section \ref{section:groupactions}.
Finally we can relate $HF_*(\phi,1)$ with $HF^*(\phi)$:
\begin{theorem} \label{thm:floerlongexactsequence}
For any symplectomorphism $\phi$, we have a long exact sequence
\[\rightarrow HF^i(\phi) \rightarrow HF_i(\phi,1) \rightarrow
HF^{i-1}(\phi) \rightarrow. \]
\end{theorem}
This is proven is section \ref{section:floerlongexactsequence}.
It is very similar to the Gysin exact sequence
from \cite{BourgeoisOancea:exactsequence}.
These theorems might give us information about symplectic homology
if we know a lot about Floer homology of symplectomorphisms.
For instance the work from \cite{CottonClay:symplecticfloerarea}
enables us to calculate these groups if $F$ is a surface.
In this paper, we go the other way around.
We use symplectic homology $SH_*(E)$ to give us information about $\phi$.

We have the following application of the previous three theorems:
For any symplectomorphism $\phi$, we can define its {\it positive stabilization}.
This involves adding a Weinstein $n$-handle to $F$ to create $F'$ and
changing $\phi$ to $\tau_S \circ \phi$ where $\tau_S$ is the Dehn
twist about an exact Lagrangian sphere which intersects the cocore
of the handle in exactly one point.
This construction is described in more detail at the start of section
\ref{section:appendixmainargument}.
Also see \cite[Section 1]{Seidel:longexactsequence} for a detailed discussion
of Dehn twists and their relation with Lefschetz fibrations and
\cite[Section 2.2]{Cieliebak:handleattach} or \cite{Weinstein:contactsurgery}
for a description of Weinstein handle attaching.
\begin{corollary} \label{corollary:nonvanishingfloerhomology}
Suppose that $\phi : F \rightarrow F$ is a symplectomorphism
such that it is
obtained by one or more stabilizations to the identity map
$\text{id} : F' \rightarrow F'$. Then if 
the Euler characteristic is odd,
then $HF^*(\phi^k,\Q) \neq 0$ for infinitely
many $k$.
\end{corollary}
This implies that
for infinitely many $k$, any symplectomorphism $\psi : F \rightarrow F$
which is the identity at the boundary and Hamiltonian isotopic through such symplectomorphisms
to $\phi^k$ has at least one fixed point away from the boundary.
This corollary is proven in section \ref{section:applications}.
The key idea here is that $\phi$ is the monodromy
map of some Lefschetz fibration symplectomoprhic to
$\C \times \widehat{F}$. This has symplectic homology zero.
We also use the fact that the $E^\infty$
page of a spectral sequence is non-trivial if
the total dimension of the $E^1$ pages is odd.

We also have the following corollary of Theorems \ref{thm:spectralsequence},
\ref{thm:groupactions} and \ref{thm:floerlongexactsequence}:
\begin{corollary} \label{thm:infiniterank}
Suppose $E$ is the total space of a Lefschetz fibration with monodromy $\phi$
such that $SH_*(E,\Q)$ has infinite rank,
then $HF^*(\phi^k,\Q)$ is nonzero for infinitely many $k$.
\end{corollary}
This is basically because $HF^*(\phi^k)$ is finite
dimensional for all $k$, but symplectic homology
is infinite dimensional. Hence $HF^*(\phi^k)$
must be non-trivial for infinitely many $k$
because we need $\bigoplus_k HF^*(\phi^k)$
to be infinite dimensional.
There are many examples of Liouville domains with infinite
dimensional symplectic homology such as cotangent
bundles of manifolds and also there
are examples from \cite{McLeanAlbers:leafwise}.
Many of these examples admit Lefschetz fibrations.
In fact all of these have boundaries admitting open book
decompositions (see \cite{Giroux:openbooks} and \cite{Akbulut:lefschetzfibrationsoncompactstein}), and we can generalize
some of the above theorems in the following way:

Suppose $p : \partial M \setminus B \twoheadrightarrow S^1$ is an open book on the contact boundary of a Liouville domain $M$.
An open book is described as follows:
The map $p$ is a fibration whose fibers are symplectomorphic
to the interiors of Liouville domains (in fact these are a special
class of Liouville domains called Stein domains).
The symplectic form on these fibers is the one induced from $d\alpha_M$
where $\alpha_M$ is the contact form on $\partial M$.
The closure of each fiber is a Liouville domain with Liouville form
$\alpha_M|_{\text{fiber}}$ and the boundary is $B$.
This has a natural connection and a monodromy map
coming from parallel transport around the base $S^1$.
This monodromy map is some symplectomorphism $\phi : F \rightarrow F$ fixing the boundary.
If we look at the statement of Theorems \ref{thm:spectralsequence} and \ref{thm:infiniterank},
then if we change $\phi$ so that it is now the monodromy map of our open book,
then they are still true as long as we have the additional condition that
the natural map $H^1(M) \rightarrow H^1(\partial M)$ is surjective.
This is because every open book decomposition of $\partial M$
gives us a map $\pi : \widehat{M} \setminus K \rightarrow \C$
where $K$ is a compact subset of $\widehat{M}$
and such that $\pi$ looks like a Lefschetz fibration away from $K$.
This map $\pi$ has a monodromy map equal to the monodromy
map $\phi$ of the open book decomposition.
The point is that when we prove Theorem \ref{thm:spectralsequence},
we do not use the interior of the Lefschetz fibration,
we only use the boundary.
Some of these issues
will be dealt with in a future paper \cite{McLean:computability}.

\bigskip

{\bf Acknowledgments:}
I would like to thank Ivan Smith, Peter Albers and Paul
Seidel for useful comments. 
The author was partially supported by
NSF grant DMS-1005365.

\section{Definitions of all our Floer homology groups}
\subsection{Lefschetz fibrations with positive and negative ends} \label{section:lefschetzfibrations}

Before we define other Floer homology groups
we will define Lefschetz fibrations with positive and negative ends.
These Lefschetz fibrations will also be used in the definition
of the Floer  homology groups  $HF_*(\phi,k)$.
Also we use these fibrations to construct our spectral sequence.

Let $\phi : \widehat{F} \rightarrow \widehat{F}$
be a compactly supported symplectomorphism.
Any symplectomorphism $\phi$ is {\it exact} if $\phi^*(\theta_F) = \theta_F + df$
where $\theta_F$ is the Liouville form and $f : \widehat{F} \rightarrow \R$
is a function.
We have that $\phi$ is isotopic to an exact symplectomorphism
(see the proof of \cite[Lemma 1.1]{BEE:legendriansurgery}).
So from now on by this deformation argument we can assume that $\phi$ is exact.
Before we define Lefschetz fibrations, we will construct
for each exact symplectomorphism $\phi$ a mapping torus $M_\phi$.
The mapping torus is a fibration $\pi_\phi : M_\phi \twoheadrightarrow S^1$
with the following properties:
\begin{enumerate}
\item There exists a contact form $\alpha_\phi$ on $M_\phi$.
\item $d\alpha_\phi$ restricted to each fiber is a symplectic
form. This means we have a connection on this fibration
coming from the line field that is $d\alpha_\phi$ orthogonal
to the fibers.
\item The monodromy map going positively around $S^1$ is Hamiltonian isotopic to $\phi$.
This means that it is equal to $\phi$ composed with a compactly supported Hamiltonian symplectomorphism.
\item \label{item:propertynearinfinity}
Near infinity, the fibration is equal to
the product fibration
\[ [R,\infty) \times \partial F \times S^1 \twoheadrightarrow S^1\]
where $\alpha = d\theta + \theta_F$.
Here $\theta$ is the angle coordinate.
\end{enumerate}
These properties determine $M_\phi$ up to isotopy of contact structures
(in fact up to contact isomorphism by Gray's stability theorem
even though the contact manifold is non-compact).

Here is an explicit construction of $M_\phi$:
As a manifold, $M_\phi$ is equal to
$[0,1] \times \widehat{F} / \sim$
where the equivalence relation $\sim$
identifies $(0,x)$ with $(1,\phi(x))$.
We have a contact form $dt + \theta_F$
on the product $[0,2\pi] \times \widehat{F}$
where $t$ parameterizes $[0,2\pi]$.
The problem is that this does not fit well with
our equivalence relation
so we cannot pass directly to the quotient.
The symplectomorphism $\phi$ is exact, which means
that $\phi^*\theta_F = \theta_F + df$.
In our case we can assume that $f=0$ outside a compact set.
Choose a function
$G : [0,2\pi] \times \widehat{F} \rightarrow \R$
such that
$G(t,x) = 0$ near $t=0$ and $G(t,x) = f(x)$ near  $t=2\pi$.
We also require that $G= 0$ outside a large compact set.
For a constant $C$ large enough,
$Cdt + \theta_F + dG$ is a contact form.
This also passes to the quotient and
hence defines a contact form $\alpha'_\phi$.
This satisfies all the properties stated above except
possibly the last one. The contact form $\alpha'_\phi$
can be rescaled by a constant so that it
satisfies property (\ref{item:propertynearinfinity})
for the following reason:
Near infinity, we have that $\alpha'_\phi$
is equal to $Cdt + \theta_F$,
where $\theta_F = r_F \alpha_F$ where $\alpha_F$
is the contact form on $\partial_F$.
If we multiply $\alpha'_\phi$ by $\frac{1}{C}$, then
it must look like $dt + \frac{1}{C} r_F \alpha_F$ near infinity.
We can construct a diffeomorphism
$\Phi : M_\phi \rightarrow M_\phi$ such that
outside the region 
\[\left\{ r_F \geq 1 \right\} = [1,\infty) \times \partial F \times S^1,\]
we have that $\Phi = \text{id}$ and
inside this region $\Phi(x,y,z) = (h(x),y,z)$ for some $h:[1,\infty) \rightarrow [1,\infty)$.
We define $h$ so that $h(r_F) = \text{id}$ near $r_F = 1$
and $h_F(r_F) = Cr_F$ for $r_F$ large.
Then $\alpha_\phi := \Phi^*(\alpha'_\phi)$ is equal to $dt + \theta_F$
near infinity. This satisfies all the required properties.

Let $S$ be an oriented surface with $s_+ + s_-$ punctures
where we label $s_+$ of these punctures with a $+$ sign
and the others with a $-$ sign. The $+$ punctures are
called {\it positive punctures} and the $-$
punctures are called {\it negative punctures}.
Around each positive puncture, we choose an orientation preserving diffeomorphism
from $[1,\infty) \times S^1$ to a neighborhood of this puncture.
Also around each negative puncture we choose an orientation preserving diffeomorphism
from $(0,1] \times S^1$ to a neighborhood of this puncture.
These diffeomorphisms are called {\it framings}
around each puncture.
On $S$ we choose a $1$-form $\theta_S$ such that
$d\theta_S$ is a symplectic form
and such that $\theta_S = rd\theta$
around each puncture where
$r$ and $\theta$ parameterizes
$[1,\infty)$ (or $(0,1]$) and $S^1$ respectively.
It turns out that for any oriented punctured surface
with at least one positive labeled puncture, we can find such a structure
and this structure is unique up to isotopy.

Let $\pi : E \twoheadrightarrow S$ be a smooth map with finitely many singularities
such that the derivative of $\pi$ is surjective  away from these singularities.
We let $E$ have a symplectic form $\omega_E$ making all the smooth fibers into symplectic manifolds
and that $\omega_E = d\theta_E$ for some $1$-form $\theta_E$.
We assume that $\pi$ has the following properties:
\begin{enumerate}
\item The smooth fibers of $E$
are symplectomorphic to completion $\widehat{F}$ of a Liouville domain $(F,\theta_F)$.
\item There exists a subset $E_h \subset E$ such that $(E \setminus E_h) \cap \pi^{-1}(p)$
is relatively compact for all $p \in S$ and such that
$\pi|_{E_h}$ is the trivial fibration
\[S \times \partial(F) \times [1,\infty) \twoheadrightarrow S.\]
Also on $E_h$, $\theta_E = \theta_S + \theta_F$.
\item Around each positive puncture, $\pi$ is equal to the product
\[(\text{id},\pi_\phi) : [1,\infty) \times M_\phi \twoheadrightarrow [1,\infty) \times S^1\]
where $\pi_\phi : M_\phi \twoheadrightarrow S^1$ is the mapping torus of some
symplectomorphism $\phi$.
We assume that $\theta_E$ is equal to $r \wedge \pi_\phi^* \vartheta + \alpha$
where $\alpha$ is a contact form on $M_\phi$ and $\vartheta$ is the angle coordinate on $S^1$.
This is called a {\it positive cylindrical end}.
\item We have a similar model $(0,1] \times M_\phi$ around each negative
puncture. This is called a {\it negative cylindrical end}.
\item
Around each singularity $p$, there exists a complex structure $J_p$ compatible
with the symplectic form, and another complex structure $j_p$ around $\pi(p)$
making $p$ holomorphic and equal to \[(z_1,\cdots,z_n) \twoheadrightarrow \sum z_i^2\]
for some chosen holomorphic charts around $p$ and $\pi(p)$.
\end{enumerate}
We say that an almost complex structures $J$ is {\it compatible}
with a symplectic form $\omega$ if:
\begin{enumerate}
\item $\omega(\cdot,J(\cdot))$ is symmetric and positive definite.
\item $\omega(J(\cdot),J(\cdot)) = \omega(\cdot,\cdot)$.
\end{enumerate}

\begin{defn} \label{defn:lefschetzfibration}
The above map $\pi : E \rightarrow S$ is called a Lefschetz
fibration with $s_-$ negative ends and $s_+$ positive ends.
An {\it isotopy of Lefschetz fibrations}
is a smooth family
of Lefschetz fibrations parameterized by $[0,1]$.
This means that all the data such as the $1$-form $\theta_E$,
almost complex structure, the choice of trivialization
of $\pi|_{E_h}$ and the cylindrical ends varies smoothly.
\end{defn}

We are interested in the Lefschetz fibrations up to isotopy.
If $p$ is a regular point of $\pi$, then the tangent space at $p$
splits up as $TE^v \bigoplus TE^h$ where $TE^v$ is the set
of vectors tangent to the fiber at $p$ and $TE^h$ is the
set of vectors $\omega_E$ orthogonal to $TE^h$.
Each Lefschetz fibration has a connection (well defined away from the singularities)
given by the plane distribution $TE^h$.
The parallel transport maps are exact symplectomorphisms.
If we are given a small circle around some puncture $p$ which winds
in the direction where $\theta$ is increasing,
then this induces a compactly supported symplectomorphism $\phi_p : F \rightarrow F$
called the monodromy around $p$.
This is unique up to isotopy.

\subsection{Definitions of Floer cohomology for symplectomorphisms} \label{section:floerhomologydefn}

We will use coefficients in a field $\K$.
In general our homology theory works if we have coefficients in a principal ideal domain.
Here we define the Floer cohomology groups $HF_*(\phi)$,
where $\phi : \widehat{F} \rightarrow \widehat{F}$ is a compactly supported symplectomorphism
and where $\widehat{F}$ is the completion of a Liouville domain $(F,\theta_F)$.
We write $\omega_F := d\theta_F$.
The boundary $\partial F$ has a natural contact form $\alpha_F := \theta_F|_{\partial F}$
and $\theta_F = r_F\alpha_F$ on the cylindrical end $[1,\infty) \times \partial F$
where $r_F$ parameterizes $[1,\infty)$.
For simplicity we assume that the first Chern class of the
symplectic manifold $F$ is trivial.
We will assume that $\phi$ is an exact symplectomorphism.
Any compactly supported symplectomorphism is isotopic
through compactly supported symplectomorphisms to an exact
symplectomorphism anyway so this does not really
put any constraint on $\phi$
(see the proof of \cite[Lemma 1.1]{BEE:legendriansurgery}).
%
Let $H : \widehat{F} \rightarrow \R$ be a Hamiltonian such that
$H|_F = 0$ and $H = h(r_F)$ on the cylindrical end of $\widehat{F}$.
On the manifold $F_R := \{r_F = R\}$, the Hamiltonian flow is equal to
$-Rh'(R)X$ where $X$ is the Reeb flow of $\partial F \cong \{r_F = R\}$
where the diffeomorphism is induced from
the projection $[1,\infty) \times \partial F \twoheadrightarrow \partial F$.
We assume that
$H = 0$ on the region where $\phi$ is non-trivial
and that $h'(r)$ is so small
that $H$ has no $1$-periodic orbits in the region where $h'(r) \neq 0$.
We also require that for $r \gg 0$, $h'(r) > 0$.
Let $\phi_H^1$ be the time $1$ Hamiltonian flow of this Hamiltonian.
This means that the symplectomorphism $\phi \circ \phi^1_H$
has all its fixed points lying inside a compact subset of $\widehat{F}$.
We say that a point $x$ is a {\it non-degenerate fixed point}
if $\phi(x) = x$ and $D_x\phi : T_xM \rightarrow T_xM$
has no eigenvalues equal to $1$.
By using work from \cite{DostoglouSalamon:instantonscurves},
we can perturb $\phi \circ \phi^1_H$ by a generic
compactly supported Hamiltonian symplectomorphism
$\phi^1_{H'}$ so that all the periodic orbits
of $\phi' := \phi \circ \phi^1_H \circ \phi^1_{H'}$
are non-degenerate.
In particular this means that there are only finitely
many $1$-periodic orbits.
\begin{defn} \label{label:standardperturbation}
Any symplectomorphism $\phi'$ constructed as above
from $\phi$ is called a standard perturbation of $\phi$.
\end{defn}

Let $C_*(\phi')$ be the free $\K$ vector space generated
by the fixed points of $\phi'$.
This is a graded vector space, and its grading
is the Conley-Zehnder index
taken with negative sign.
Choose an almost complex structure $J$ on $T\widehat{F}$.
We have a complex line bundle $\Lambda^n(T\widehat{F},J)$ given by the highest exterior
power of the complex vector bundle $(T\widehat{F},J)$.

In order to define the Conley-Zehnder index for fixed points
of a symplectomorphism, we need to choose
a smooth family of bundle trivializations as follows:
Let $J_t$ be a family of almost complex structures
compatible with the symplectic form $\omega_M$ on $M$
parameterized by $t \in [0,1]$.
We also assume that $J_1 \circ D\phi' = D\phi' \circ J_0$
(i.e. $\phi'$ is $(J_0,J_1)$-holomorphic).
These almost complex structures
$J_t$ must by {\it cylindrical at infinity}
which means that outside a large compact set,
$dr_F \circ J_t = -\theta_F$.
We also assume that $J_t$ is independent
of $t$ near infinity.
Our series of bundle trivializations are
\[\text{Tr}(t) : \widehat{F} \times \C \rightarrow \Lambda^n(T\widehat{F},J_t)\]
parameterized by $t \in [0,1]$.
Because $\phi'$ is $(J_0,J_1)$-holomorphic, we have that $D\phi'$
induces a map $\Lambda^n D\phi'$ from $\Lambda^n(T\widehat{F},J_0)$ to
$\Lambda^n(T\widehat{F},J_1)$ given by
\[\Lambda^n D\phi' (v_1 \wedge v_2 \cdots \wedge v_n) := D\phi'(v_1) \wedge \cdots \wedge D\phi'(v_n).\]
We assume that the sequence of trivializations
above fits into the following commutative diagram:
\[
\xy
(0,0)*{}="A"; (40,0)*{}="B";
(0,-20)*{}="C"; (40,-20)*{}="D";
"A" *{\widehat{F} \times \C};
"B" *{\Lambda^n(T\widehat{F},J_0)};
"C" *{\widehat{F} \times \C};
"D" *{\Lambda^n(T\widehat{F},J_1)};
{\ar@{->} "A"+(10,0)*{};"B"-(10,0)*{}};
{\ar@{->} "C"+(10,0)*{};"D"-(10,0)*{}};
{\ar@{->} "A"+(0,-4)*{};"C"+(0,4)*{}};
{\ar@{->} "B"+(0,-4)*{};"D"+(0,4)*{}};
"A"+(20,3) *{\text{Tr}(0)};
"C"+(20,3) *{\text{Tr}(1)};
"A"+(7,-10) *{\phi' \times \text{id}};
"B"+(7,-10) *{\Lambda^n D\phi'};
\endxy
\]

Another way of thinking about this family of trivializations is as follows:
The family of almost complex structure $J_t$ in fact induces an
almost complex structure on the contact distribution of the mapping
torus $M_{\phi'}$. The above set of trivializations is equivalent to a
trivialization of the highest complex exterior power of this contact distribution.

Let $x$ be a fixed point of $\phi'$. We wish to assign an index to this
point.
The trivializations  $\text{Tr}(t)$ give us a sequence of maps
\[\text{Tr}(t)_x : \C \rightarrow \Lambda^n(T_x \widehat{F},J_t).\]
We can find smooth family of complex linear maps
\[a(t) : \C^n \rightarrow T_x (\widehat{F},J_t)\]
which are also symplectomorphisms where $\C^n$ has the
standard symplectic structure such that
\begin{enumerate}
\item
the induced map
\[\Lambda^n a(t) : \C \rightarrow \Lambda^n T_x (\widehat{F},J_t)\]
is equal to $\text{Tr}(t)_x$
\item
$a(1) \circ a(0)^{-1} = D\phi'$.
\end{enumerate}
Such a choice is unique up to homotopy relative to the endpoints.
The maps $a(0)^{-1} \circ a(t)$
give us a path of symplectic matrices and then we can
use the work of \cite{RobbinSalamon:maslov} to compute an index.
For a fixed point $x$, we write $\text{ind}(x)$ for its index.

The index makes $C_*(\phi')$ into a graded vector space.
The non-degeneracy assumption combined with the fact that
all the fixed points are inside some compact set
ensures that this is a finite dimensional vector space.
We need to define a differential on this vector space.
For a fixed point $x$ of index $i$ we define $\partial x$ as follows:
\[\partial x = \sum_{\text{ind}(y)=i-1} \# \left( \M(x,y)/\R \right) y.\]
We will now define the manifold $\M(x,y)$ together
with its free $\R$ action.

The set $\M(x,y)$ consists of maps
$u : \R \times [0,1] \rightarrow M$
which are $J_t$ holomorphic such that
$u(s,t)$ tends to $x$ as $s$ tends to $-\infty$
and $y$ as $s$ tends to $+\infty$.
The map $u$ is $J_t$ holomorphic if it satisfies the following
equation:
\[u_s + J_t u_t = 0\]
where $(s,t)$ parameterizes $\R \times [0,1]$.
We also require that $u(1,t) = \phi(u(0,t))$.
The set $\M(x,y)$ has an $\R$ action given by translation
in the $s$ variable.
We have that this set $\M(x,y)/\R$ is a finite set
for generic $J_t$.
This space $\M(x,y) / \R$ is a $0$ dimensional compact oriented manifold.
Each point comes with a positive or negative orientation.
We define $\#\M(x,y)/\R$ to be equal to the number of positively
oriented points minus the number of negatively oriented points.
\begin{theorem}
\cite{DostoglouSalamon:instantonscurves}
We have $\partial^2 = 0$.
\end{theorem}
This group is independent of all choices made
(except the choice of trivializations in order to
defined the index, but we will suppress this choice here).
This gives us our Floer cohomology group $HF^*(\phi)$
for a symplectomorphism $\phi$.
%
In the introduction, our symplectomorphism was not a compactly supported symplectomorphism
in $\widehat{F}$, instead it was a symplectomorphism
from $F$ to $F$ fixing its boundary.
Because $\phi$ fixes the boundary, we have that $D\phi : TF|_{\partial F} \rightarrow TF|_{\partial F}$
is the identity. 
This is because there is a collar neighborhood
$(1-\epsilon,1] \times \partial F$
of $\partial F$ such that $\lambda = r\alpha$ where
$r$ parameterizes $(1-\epsilon,1]$, and $D\phi|_{\partial F}$ preserves
the Reeb flow of $\alpha$. Hence it must also preserve
$\frac{\partial}{\partial r}$ which is transverse to the boundary.
This implies that we can extend $\phi$ to a $C^1$ symplectomorphism
$\widehat{\phi} : \widehat{F} \rightarrow \widehat{F}$ which is compactly supported.
Hence we can define Floer homology for $\phi$ (maybe after perturbing it slightly so that
it becomes $C^\infty$ without adding extra orbits).
Hence we define $HF^*(\phi) := HF^*(\widehat{\phi})$.

\subsection{Defining another Floer homology group} \label{section:equivariantfloer}

We will first describe a family of Hamiltonians which
behave well with respect to the Lefschetz fibration $E$.
Let $M_\phi$ be the mapping torus of $\phi$.
The Reeb orbits of the contact form $\alpha_\phi$
are in one to one correspondence with fixed points
of $\phi^k$ for each $k \in \N_{>0}$.
Let $R > 0$ be a constant such that
in the region $\{r_F \geq R\}$,
we have that $\phi$ is the constant symplectomorphism.
We denote ${\mathcal R}_{\alpha_\phi}$ to be
the set of Reeb orbits of $\alpha_\phi$
corresponding to fixed points of $\phi$
outside the region $\{r_F \geq R\}$.
By using a similar argument to \cite{} can perturb the
contact form $\alpha_\phi$ generically inside the
region $\{r_F < R\}$
so that all the Reeb orbits inside this region are non-degenerate
(the point here is that we consider the moduli space of
contact forms $\alpha$ such that $\alpha = \alpha_\phi$
in the region $\{r_F \geq R\}$).
The {\it period spectrum} is a subset of $\R$ corresponding
to the set of lengths
of Reeb orbits.

Let $E$ be a Lefschetz fibration as described in section
\ref{section:lefschetzfibrations}.
Each cylindrical end is either of the form
$[1,\infty) \times M_\phi$ or
$(0,1] \times M_\phi$.


For each function $g : (0,\infty) \rightarrow \R$
such that $g(r) = 0$ for $r$ near $1$,
we will construct a Hamiltonian $H_g : E \rightarrow \R$ as follows:
Let $S$ be the base of the Hamiltonian.
The region $E_h$ as described in section \ref{section:lefschetzfibrations}
is equal to $S \times [1,\infty) \times \partial F$.
We choose a Hamiltonian $H : E_h \rightarrow \R$
such that $H = h(r_F)$ where $r_F$ parameterizes
the interval $[1,\infty)$.
We choose $h$ so that $h(r_F) = 0$ near $r_F = 1$,
and $h'(r_F) > 0$ near infinity.
We also ensure that $h'(r_F)$ is small so that in the region
where $h'(r_F) \neq 0$, $H$ has no periodic orbits.
We also require $h<1$.
The Hamiltonian $H$ can be extended to a Hamiltonian
$\widetilde{H} : E \rightarrow \R$ which is equal to $0$
outside $E_h$. This is because $H = 0$ for near $r_F = 1$.

Let $I$ be the interval $[1,\infty)$ or $(0,1]$.
Let $C$ be a cylindrical end of $E$.
This cylindrical end is of the form $I \times M_\phi$.
We define $H^C_g$ as $g(s) + \widetilde{H}$ on this cylindrical end
for any $\widetilde{H}$ described as above.
We define $H_g$ to be equal to
$\widetilde{H}$ outside all the cylindrical ends,
and for each cylindrical end $C$, we define
$H_g$ to be equal to $H_g^C$.
This function is smooth because $g(r) = 0$
for $r$ close to $1$.

From now on we assume that $g'(r)$ is constant and not a multiple of $2\pi$ when $r \gg 0$.
On any cylindrical end, the Liouville form is equal to $r d\vartheta + \alpha$
hence the Hamiltonian vector field of $g(r)$ is equal to
the horizontal lift of the Hamiltonian vector field of $g(r)$
on $[1,\infty) \times S^1$ with symplectic form $dr \wedge d\alpha$.
Hence $g(r)$ has no $1$-periodic orbits for $r \gg 0$ because $g'(r)$ is not
a multiple of $2\pi$.
Also, we assume that $g'(r) > 0$ and very small (less than $2\pi$ will do) for
$r$ very small and also that $g'$ is constant on some open subset of this region.
This ensures that $H_g$ has no
$1$-periodic orbits on the negative cylindrical ends when $r$ is small.
\begin{defn} \label{defn:lefschetzadmissiblehamiltonian}
A Hamiltonian $H : S^1 \times E \rightarrow \R$ is
{\it Lefschetz admissible of slope} $\lambda$
if there exists a $g : (0,\infty) \rightarrow \R$
as described above such that
$H = H_g$ outside some compact set $K$.
\end{defn}
A Hamiltonian $H$ is said to be non-degenerate
if all the fixed points of the symplectomorphism
$\phi^1_H$ are non-degenerate.
The set of non-degenerate Lefschetz admissible Hamiltonians
is generic among Lefschetz admissible Hamiltonians.

The map $\phi^1_{H_g}$ is a symplectomorphism
so we can define $HF_*(\phi^1_{H_g})$ as above.
The only difference is that the symplectomorphism $\phi^1_{H_g}$
is not compactly supported, but this does not affect the definition.
For a Hamiltonian symplectomorphism, we can define $HF_*(\phi^1_{H_g})$
in a slightly different way. This definition will be
useful for our purposes even though the definition is almost identical.
There is a $1$-$1$ correspondence between fixed points
$x$ of $\phi^1_{H_g}$ and $1$ periodic orbits $\bar{x} : S^1 \rightarrow E$ of
the Hamiltonian flow $X_H$ (I.e. the vector field $X_H$ satisfies $\omega(X_H,\cdot) = dH$).
Hence the chain complex $C_*(\phi^1_H)$
is generated by $1$-periodic orbits of $X_H$.
If we choose a trivialization of the line bundle $\Lambda^{n+1} TE$
(where the dimension of $E$ is $2n +2$), then we have a well defined grading.
We choose the trivialization as follows:
Let $j$ be the natural complex structure on $\R \times S^1$
where $j(\frac{\partial}{\partial s}) = \frac{\partial}{\partial t}$.
Let $J$ be an almost complex structure compatible
with the symplectic form $\omega_E$
making $\pi_E$ holomorphic.
The vertical tangent spaces
$TE^v$ are complex vector bundles because $\pi$
is $(J,j)$-holomorphic.
The tangent bundle of the base $\R \times S^1$ of the Lefschetz fibration
also has a natural trivialization, and hence
the horizontal plane bundle $TE^h$ has a natural trivialization.
Combining these two trivializations gives us a trivialization
of $\Lambda^n TE^v \bigoplus TE^h$ and hence
we get a trivialization of $\Lambda^{n+1} TE$.
This trivialization enables us to define
for each $1$-periodic orbit $\bar{x}$ an
index which we denote by $\text{ind}(\bar{x})$.

We now need to define a differential.
In order to do this, we first need to describe
the set of Lefschetz admissible complex structures
$\J(\pi_E)$.
Let $J$ be an almost complex structure making
$\pi_E$ holomorphic.
On the region $E_h = S \times [1,\infty) \times \partial F$,
we let $J = j \bigoplus J_F$ where
$J_F$ is some almost complex structure on $\widehat{F}$
that is cylindrical at infinity.
We also assume that on each cylindrical end $I \times M_\phi$,
that $J$ is invariant under translations in the $r$ direction
where $r$ parameterizes $I$.
\begin{defn} \label{defn:lefschetzadmissiblealmostcomplexstructure}
An $S^1$ family of almost complex structures $J_t$ compatible with $\omega_E$
is {\it Lefschetz admissible} if outside some compact set $K$,
$J_t = J$ for some $J$ as described above.
\end{defn}

The differential is a count of cylinders satisfying the perturbed
Cauchy Riemann equations.
We define the differential $\partial$ as follows:
\[\partial \bar{x} = \sum_{\text{ind}(\bar{y})=\text{ind}(\bar{x})-1}
\left(\#\M(\bar{x},\bar{y})/\R\right) \bar{y}.\]
Here $\M(\bar{x},\bar{y})$
is the set of maps $u : \R \times S^1 \rightarrow E$
satisfying
\[\partial_s u + J_t (\partial_t u - X_H) = 0\]
where $X_H$ is the Hamiltonian vector field of $H$.
We also require that $u(s,t)$ tends to $\bar{x}(t)$
as $s$ tends to $-\infty$ and to $\bar{y}(t)$
as $s$ tends to $+\infty$.
The $\R$ action again is translations in the $s$ coordinate.
In order for $\partial$ to be well defined and a differential,
we need to ensure that these moduli spaces are compact.
If all these solutions stay inside a compact subset of $E$,
then the results from \cite{BEHWZ:compactnessfieldtheory} ensure
that $\#\M(\bar{x},\bar{y})/\R$ is compact.
The problem is that $E$ is non-compact and these
solutions $u$ could escape to infinity.
The maximum principle \cite[Lemma 7.2]{SeidelAbouzaid:viterbo}
combined with \cite[Lemma 5.2]{McLean:symhomlef}
and Lemma \ref{lemma:lefschetzminimumprinciple}
ensures that this does not happen.

For generic $(H_t,J_t)$, the map
$\partial$ is well defined
and satisfies $\partial^2 = 0$.
If this is true then we say that
the pair $(H_t,J_t)$ is regular
(actually we say that this pair is regular
if a particular linear map which won't be described
here is surjective).
The homology of this chain complex is denoted by
$HF_*(\phi^1_{H_t})$.
For each element $\alpha \in H_1(E)$,
we have a subgroup $HF^\alpha_*(\phi^1_{H_t})$
generated by $1$-periodic orbits in the homology class
$\alpha$.

If we have two regular pairs $(H^1_t,J^1_t)$,
$(H^2_t,J^2_t)$
such that $H^1_t \leq H^2_t$,
then there is a natural map from
$HF_*(\phi^1_{H^1_t})$
to $HF_*(\phi^1_{H^2_t})$.
Again this also induces a map from
$HF^\alpha_*(\phi^1_{H^1_t})$
to $HF^\alpha_*(\phi^1_{H^2_t})$.
This is called the continuation map.
The idea is that we choose an increasing
homotopy $K^s_t$ of Lefschetz admissible Hamiltonians such that $K^s_t = H^1_t$
for $s \ll 0$ and $K^s_t = H^2_t$ for $s \gg 0$.
We also choose a smooth family of Lefschetz admissible
almost complex structures ${J'}^s_t$ such that
${J'}^s_t = J^1_t$ for $t \ll 0$ and
${J'}^s_t = J^2_t$ for $t \gg 0$.
The continuation map sends an orbit $\bar{x}$ of $H^1$
to:
\[\sum_{\text{ind}(\bar{y}) = \text{ind}(\bar{x})}
\#\M(\bar{x},\bar{y},{J'}^s_t) \bar{y}.\]
The space $\M(\bar{x},\bar{y},{J'}^s_t)$
is the set of maps $u : S^1 \times \R \rightarrow E$
satisfying
\[\partial_s u + {J'}^s_t (\partial_t u - X_{K^s_t}) = 0\]
joining $\bar{x}$ and $\bar{y}$.
This map is well defined only if we have an
non-decreasing homotopy.
The maximum principles
\cite{SeidelAbouzaid:viterbo}
combined with \cite[Lemma 5.2]{McLean:symhomlef}
and Lemma \ref{lemma:lefschetzminimumprinciple}.
If we have a different family of Hamiltonians $K^s_t$
and almost complex structures ${J'}^s_t$ then
we get a chain homotopic map.
This means that we get the same map on homology.
If we compose two continuation maps, then we get
a continuation map.

For each orbit, we have its action which is defined
as
\[{\mathcal A}(x) := -\int_x \theta_E - \int_x H dt.\]
The differential $\partial$ sends an orbit
of action $A$ to a linear combination of orbits
of action $\leq A$.
Hence we can define a subcomplex
$C_*^{(-\infty,a)}(H)$
generated by orbits
of action strictly less than $a$.
We also have a quotient complex:
\[C_*^{[b,a)}(H) := C_*^{(-\infty,b)}(H) / C_*^{(-\infty,a)}(H).\]
The homology of this complex is denoted by
$HF_*^{[b,a)}(H)$.
The continuation maps also decrease action,
and hence they induce chain maps:
\[C_*^{[b,a)}(H_1) \rightarrow C_*^{[b,a)}(H_2).\]
For a Lefschetz fibration $E$ and a subset $K \subset E$,
we have that the set of regular Lefschetz admissible pairs
$(H,J)$ form a directed system with respect
to the ordering $\leq$ where $(H_1,J_1) \leq (H_2,J_2)$
if and only if $H_1 \leq H_2$.
This means that the groups
$HF_*^{[b,a)}(H,J)$ also form a directed system
where the morphisms are continuation maps.
We define
\[SH_*^+(\pi_E,K) := \varinjlim_{(H,J),H|_K < 0} HF_*^{[0,\infty)}(H,J).\]
The direct limit is taken over Lefschetz admissible pairs $(H,J)$.
We also define \[SH_*(\pi_E) := \varinjlim_{(H,J)} HF_*(H,J).\]
Really these direct limits should be done on the chain
level, but this produces the same result.
If $\alpha \subset H_1(E)$, we can also define
$SH_*^\alpha(E,K)$ by only considering orbits
in the $H_1$ class is contained in $\alpha$.

If $K_1 \subset K_2$, we have a natural morphism
$SH_*(\pi_E,K_2) \rightarrow SH_*(\pi_E,K_1)$.
This is called a {\it transfer map}.
This is because the directed system of pairs $(H,J)$
defining $SH_*(\pi_E,K_2)$
is a subdirected system of the one defining
$SH_*(\pi_E,K_1)$.

Let $\phi$ be a symplectomorphism.
Let $W_\phi$ be the Lefschetz fibration
$(0,\infty) \times M_\phi$ with $1$-form
$\theta_{\phi} := s dp_\phi^* \vartheta + \alpha_\phi$
where $\alpha_\phi$ is the natural contact
form on $M_\phi$ and $s$ parameterizes $(0,\infty)$.
Also $\vartheta$ is the angle coordinate on $S^1$
and $p : M_\phi \twoheadrightarrow S^1$
is the natural mapping torus fibration map.
The Lefschetz fibration map
\[\pi_\phi : W_\phi = \R \times M_\phi \twoheadrightarrow \R \times S^1\]
is given by
$\pi_\phi(s,x) = (s,p(x))$. For each $k \in \N$
we define $\beta_k \subset H_1(W_\phi,\Z)$
to be the subset
$\phi_*^{-1}(\{l\})$ where $l \in H_1(\R \times S^1)$
is represented by the loop winding around the $S^1$ factor 
negatively $k$ times around (i.e. so that integrating $d\vartheta$ over
the loop is negative).
\begin{defn} \label{defn:otherfloerhomologygroup}
\[HF_*(\phi,k) := SH_*^{\beta_k}(\pi_{\phi}).\]
\end{defn}
Here is an invariance result.
\begin{lemma} \label{lemma:symplecticinvariance}
If $\phi'$
is another exact symplectomorphism isotopic to $\phi$
through compactly supported exact symplectomorphisms,
then $HF_*(\phi,k)$ is isomorphic to $HF_*(\phi',k)$.
\end{lemma}
We need some preliminary Lemmas first:
\begin{lemma} \label{lemma:mappingtorusinvariance}
We have that $(M_\phi,\alpha_\phi)$ is contactomorphic
to $(M_{\phi'},\alpha_{\phi'})$.
In fact there exists a function $f : M_\phi \rightarrow (0,\infty)$
such that $f=1$ outside a large compact set
and a diffeomorphism $P : M_{\phi} \rightarrow M_{\phi'}$ such that
$P^* \alpha_{\phi'} = f \alpha_\phi$.
\end{lemma}
\proof of Lemma \ref{lemma:mappingtorusinvariance}.
Let $\phi_t : F \rightarrow F$ be a path
parameterized by $t \in [0,1]$ of compactly supported exact symplectomorphisms
joining $\phi$ and $\phi'$.
We have that $M_{\phi}$ is diffeomorphic to $M_{\phi'}$.
So from now on we regard the contact form $\alpha_{\phi'}$
to be a contact form on $M_{\phi}$ making
$(M_{\phi},\alpha_{\phi'})$ contactomorphic to $M_{\phi'}$.
We can choose a series of contact forms $\alpha_{\phi_t}$
on $M_\phi$ such that $\alpha_{\phi_0} = \alpha_\phi$
and $\alpha_{\phi_1} = \alpha_{\phi'}$.
We can also ensure that there exists a constant $R>0$
such that the region $\left\{ r_F \geq R \right\}$
with contact form $\alpha_{\phi_t}$ is contactomorphic
to $[R,\infty) \times \partial F \times S^1$
with contact form $d\theta + r_F \theta_F$.
Hence all the contact forms are equal near infinity,
so we can use Grays stability theorem to show that
they are in fact contactomorphic.
\qed

\begin{lemma} \label{lemma:alternativefloerdefinition}
Let $H: S^1 \times W_\phi \rightarrow \R$
be a Lefschetz admissible Hamiltonian of slope greater than $2\pi k$
then
\[HF_*(\phi,k) := HF_*^{\beta_k}(H,J).\]
\end{lemma}
Because $H$ has a fixed slope, we have that
$HF_*(\phi,k)$ is a finite dimensional vector space.
\proof of Lemma \ref{lemma:alternativefloerdefinition}.
First of all note that if we have two Hamiltonians
$H_1$ and $H_2$ of the same slope then
$HF_*^{\beta_k}(H_1,J) \cong HF_*^{\beta_k}(H_2,J)$.
This is because there exists a constant $C>0$ such that
\[H_1 - C < H_2 < H_1 + C < H_2 + 2C\]
and hence we have natural continuation maps:
\[HF_*^{\beta_k}(H_1 - C,J) \rightarrow HF_*^{\beta_k}(H_2,J)
\rightarrow \]
\[HF_*^{\beta_k}(H_1 + C,J) \rightarrow HF_*^{\beta_k}(H_2 + 2C,J).\]
Composing any two of these continuation maps gives us an isomorphism
because adding a constant to a Hamiltonian does not change
the Floer equations or the orbits.
This means we get isomorphisms \[HF_*^{\beta_k}(H_2,J)
\cong HF_*^{\beta_k}(H_1 + C,J) \cong HF_*^{\beta_k}(H_1,J).\]
Choose $K > 2k\pi$ which is not a multiple of $2\pi$.
Let $h(r_F)$ be the function as defined at the start of this section.
We assume that the derivatives of $h$ are so small that
the only $1$-periodic orbits are the constant
orbits when $h=0$.

Let $g_s : (0,\infty) \rightarrow \R$
be a family of functions parameterized by $s \in [0,\infty)$ with the following properties:
\begin{enumerate}
\item $g'_s(r) < \pi$ for $r \leq 1$.
\item $g'_s>0,g_s''(s) \geq 0$.
\item $g'_s(r) = K$ for $1 \leq r \leq 2$.
\item $g'_s(r) = K+s$ for $s \geq 3$.
\item $g_s(r)$ tends to infinity pointwise, and $g'_s$ tends to $0$ in the
region $r \leq 1$ as $s$ tends to infinity.
\item $g_s(r) = g_0(r)$ for $r \leq 2$.
\end{enumerate}.
Let $J$ be an almost complex structure which is Lefschetz admissible
and makes $\pi_\phi$ $(J,j)$ holomorphic where
$j$ is the complex structure on $(0,\infty) \times S^1$
given by quotienting the upper half plane in $\C$ by $\Z$.
By definition, we get that
\[HF_*(\phi,k) = \varinjlim_s HF_*^{\beta_k}(g_s(r) + h(r_F),J).\]

Let $\widetilde{\frac{\partial}{\partial \vartheta}}$ be the horizontal
lift of the vector field $\frac{\partial}{\partial \vartheta}$
where $\vartheta$ is the angle coordinate for $S^1$.
The Hamiltonian $g_s(r) + h(r_F)$ has no $1$-periodic
orbits in the region $\{r_F \geq 1\}$
also all the orbits that wrap around the base
$k$ times are situated in the region where $g_s'(r) = 2k\pi$.
This is a subset of the region $r \leq 2$.
Hence all the generators of $HF_*^{\beta_k}(g_s(r) + h(r_F))$
are contained in this region.
The maximum principle \cite[Lemma 5.2]{McLean:symhomlef}
ensures that any Floer trajectory or continuation map trajectory
connecting orbits of $g_{s_1} + h(r_F)$ and $g_{s_2} + h(r_F)$ for $s_1 \leq s_2$
is contained in the region $r \leq 2$.
We have that $g_s = g_0$ in this region and hence
we get that $HF_*^{\beta_k}(g_s + h(r_F),J) \cong HF_*^{\beta_k}(g_0 + h(r_F),J)$
for all $J$ and all the continuation maps
\[HF_*^{\beta_k}(g_{s_1} + h(r_F),J) \rightarrow HF_*^{\beta_k}(g_{s_2} + h(r_F),J)\]
commute with this isomorphism.
Hence
\[HF_*^{\beta_k}(g_0 + h(r_F),J) = \varinjlim_s HF_*^{\beta_k}(g_s(r) + h(r_F),J)\]
and so
\[HF_*^{\beta_k}(g_0 + h(r_F),J) \cong HF_*(\phi,k).\]
This Hamiltonian is of slope $K$
and hence we have proven the Lemma.
\qed

\proof of Lemma \ref{lemma:symplecticinvariance}.
Because $M_\phi$ is contactomorphic to $M_{\phi'}$,
there is a fibrewise diffeomorphism $\Phi$ from
$W_\phi$ to $W_{\phi'}$ such that
$\Phi^*(rd\vartheta + \alpha_\phi) = rd\vartheta + f\alpha_{\phi'}$
where $f : M_\phi \rightarrow (0,\infty)$.
So from now on we will assume that $\pi_\phi$
and $\pi_{\phi'}$ are identical smooth fibrations
and $\theta_{\phi'} = r\vartheta + f\alpha_{\phi'}$.
Let $\frac{\partial'}{\partial \theta}$ be the horizontal lift of
$\frac{\partial}{\partial \theta}$ to $M_\phi$.
The function
$\alpha_\phi(\frac{\partial'}{\partial \theta})$ is positive and bounded so
has image inside $[0,Q']$ for some $Q'>0$.
Note for any function $F : W_\phi \rightarrow (0,\infty)$ such that
$F=1$ outside some compact set and such that
$\frac{\partial F}{\partial r} > -\frac{1}{Q'}$ we have that
$\theta_F := rd\vartheta + F\alpha_\phi$
is a $1$-form where
\[d\theta_F = dr \wedge \left( d\vartheta + \frac{\partial F}{\partial r} \alpha_\phi \right)
+ d(F|_{M_\phi} \alpha_\phi)\]
is a symplectic form.
Let $\widetilde{\frac{\partial}{\partial \vartheta}}$ be the horizontal
lift of the vector field $\frac{\partial}{\partial \vartheta}$
with respect to the symplectic form $d\theta_F$.
If we have a Hamiltonian of the form $h(r)$
then its Hamiltonian vector field $X^{d\theta_F}_{h(r)}$
(with respect to the symplectic form $d\theta_F$)
is 
\[-\frac{h'(r)}{1 + \frac{\partial F}{\partial r}
\alpha_\phi(\widetilde{\frac{\partial}{\partial \vartheta}}) )}
\widetilde{\frac{\partial}{\partial \vartheta}}.\]
From now on we assume that $h(r)$ is Lefschetz admissible of slope $K$.

By Lemma \ref{lemma:alternativefloerdefinition}
there exists a $K \gg 0$
(not multiple of $2\pi$) such that
for any Lefschetz admissible Hamiltonian $H$ of slope $K$,
$HF_*^{\beta_k}(H,\theta_\phi) = HF_*(\phi,k)$
and $HF_*^{\beta_k}(H,\theta_{\phi'}) = HF_*(\phi',k)$.
The function
$\alpha_\phi(\widetilde{\frac{\partial}{\partial \vartheta}})$
is positive and bounded so there is a constant
$Q>0$ such that the image of the function
is contained in $[0,Q]$.
Let $\epsilon>0$ be smaller than the smallest distance between
$K$ and any multiple of $2\pi$.
If $|\frac{\partial F}{\partial r}| < \delta := \frac{\epsilon}{Q(K - \epsilon))}$
then the function 
\[\frac{K}{1 + \frac{\partial F}{\partial r}\alpha_\phi(\widetilde{\frac{\partial}{\partial \vartheta}})}\]
is never a multiple  of $2\pi$.
In particular the Hamiltonian $Kr$ has no $1$-periodic
orbits with respect to the symplectic form $d\theta_F$.

Let $G_s : W_\phi \rightarrow (0,\infty)$ be a smooth family of functions
parameterized by $s \in [0,1]$ with the following
properties:
\begin{enumerate}
\item $\frac{\partial G_s}{\partial r} < \frac{\delta}{2}$.
\item In the region $r \leq 1$ or $r \gg 0$, $G_s = 1$.
\item $G_s = 1$ outside a compact set.
\item There is a constant $\Delta > 1$ such that
$G_1 = f$ in the region $[\Delta,\Delta+1]$.
\item $G_0 = 1$.
\end{enumerate}

Let $h : (0,\infty) \rightarrow \R$ be a function
such that
\begin{enumerate}
\item $h'(r)>0,h''(r) \geq 0$.
\item $h'(r) < \pi$ in the region $\{r \leq \Delta + \frac{1}{3}\}$
and $h'$ is constant near $\{r= \Delta\}$.
\item $h'(r) = K$ for $r \geq \Delta + \frac{2}{3}$.
\end{enumerate}
Let $J^s_t$ be a family of almost complex structures
on $W_\phi$ parameterized by $(s,t) \in [0,1] \times S^1$
such that $J^s_t$ makes $\pi_\phi$ a $(J,j)$ holomorphic
map where $j$ is a standard complex structure on the base
$(0,\infty) \times S^1 = \{\text{im}(z) > 0\} / \Z$.
We also assume that $J^s_t$ is Lefschetz admissible
with respect to the symplectic form $d\theta_{G_t}$
and the fibration $\pi_\phi$.
All the $1$-periodic orbits of the Hamiltonian $h(r)$
with respect to the symplectic form $d\theta_{G_t}$
are contained in $\{\Delta < r < \Delta +1\}$.
The maximum principle \cite[Lemma 5.2]{McLean:symhomlef}
and Lemma \ref{lemma:lefschetzminimumprinciple}
ensures that any Floer trajectory connecting these orbits
is contained in this region as well.
Also the orbits and Floer trajectories of $h(r)$
are contained in $\{\Delta < r < \Delta +1\}$
if we use the Liouville form
$\theta_f := rd\vartheta + f\alpha_\phi$.
This implies that $HF_*^{\beta_k}(h(r),\theta_{G_1})$
is isomorphic to $HF_*^{\beta_k}(h(r),\theta_f)$
which in turn is equal to $HF_*(\phi',k)$.
Also because the family of $1$-forms
$\theta_{G_t}$ are all equal to each other outside a compact set,
we get by a Moser theorem a compactly supported exact
symplectomorphism $\Phi$ from $(W_\phi,\theta_\phi)$
to $(W_\phi,\theta_{G_1})$ (which is smoothly isotopic
to the identity map).
Hence a continuation argument tells us that
\[HF_*(\phi,k) = HF_*^{\beta_k}(h(r),\theta_\phi) \cong
HF_*^{\beta_k}(\Phi^* h(r),\Phi^* \theta_{G_1})\]
\[ \cong HF_*^{\beta_k}(h(r),\theta_{G_1}) = HF_*(\phi',k).\]
This proves the Lemma.
\qed

\subsection{Definition of symplectic homology} \label{section:symplectichomologydefn}

Let $M$ be a Liouville domain with Liouville
form $\theta_M$, and let
$\widehat{M}$ be the completion of $M$.
This has a cylindrical end
diffeomorphic to $[1,\infty) \times \partial M$.
We write $r_M$  for the coordinate parameterizing
$[1,\infty)$.
We write $\{r_F \leq R\}$ for the set:
\[M \cup \left( [1,R] \times \partial M \right) \subset \widehat{M}.\]
The manifold $\partial M$ is a contact manifold
with contact form $\alpha_M := \theta_M|_{\partial M}$.
We can perturb $\theta_M$ slightly so that
the period spectrum of the contact form $\alpha_M$
is discrete and injective.

\begin{defn} \label{admissiblehamiltonian}
We say that a Hamiltonian $H : \widehat{M} \rightarrow \R$
is {\it admissible} of slope $\lambda$ if $H = \lambda r_F$
near infinity.
\end{defn}

\begin{defn} \label{admissiblealmostcomplexstructure}
An $S^1$ family almost complex structures $J_t$ compatible with
$d\theta_M$ is said to be {\it admissible}, if
$dr_F \circ J_t = -\theta_M$.
\end{defn}

All admissible pairs $(H,J)$ form a directed
system where $(H_1,J_1) \leq (H_2,J_2)$
if and only if $H_1 \leq H_2$.
We define symplectic homology as:
\[SH_*(\widehat{M}) := \varinjlim_H HF_*(H,J).\]
This is a symplectic invariant of
$\widehat{M}$.
It is invariant under exact symplectomorphisms
by \cite[section 7b]{Seidel:biasedview}.
This implies that it is in fact invariant
under general symplectomorphisms
as every symplectomorphism is isotopic
to an exact one by
\cite[Lemma 1]{BEE:legendriansurgery}.

Let $E$ be a Lefschetz fibration with one positive
cylindrical end and no negative cylindrical ends.
This is naturally symplectomorphic to the
completion of a Liouville domain as explained
in \cite[Section 1]{McLean:symhomlef}.
Let $C = [1,\infty) \times M_\phi \subset E$ be this positive
cylindrical end.
The following theorem relates symplectic homology
to Lefschetz fibrations:
Suppose that $\widehat{M}$ is symplectomorphic
to $E$.

\begin{theorem}
There is an isomorphism of groups:
\[SH_*(\widehat{M}) = SH_*(\pi).\]
\end{theorem}
This theorem is proven in \cite[Theorem 2.24]{McLean:symhomlef},
where the right hand group is called Lefschetz symplectic homology.

\section{Relationship between each of the Floer homology groups}

\subsection{Construction of the main spectral sequence} \label{section:mainspectralsequence}

In this section we prove Theorem \ref{thm:spectralsequence}.
Here is a statement of the main theorem of this section:
Let $E$ be a Lefschetz fibration with one positive
cylindrical end and no negative cylindrical ends.
This is symplectomorphic to the completion $\widehat{M}$
of a Liouville domain $M$.
Let $M_\phi$ be the mapping torus associated to the cylindrical end.
There is a spectral sequence converging to
$SH_*(\widehat{M})$ with $E^1$ page given by:
\[E_{0,q} = H^{n-*}(E)\]
And for $p > 1$,
\[E_{p,q} = HF_{q-p}(\phi,q).\]
For $p<0$, these pages are $0$.

Here is the idea of the proof of theorem \ref{thm:spectralsequence}.
We carefully construct a cofinal family of
Lefschetz admissible Hamiltonians for $SH_*(E,E \setminus C)$
where $C$ is the cylindrical end for $E$.
We ensure that
these Hamiltonians $H_i$ fix some smooth fiber $\pi^{-1}(q)$.
The base of the Lefschetz fibration is the plane $\C$.
If we project these orbits to the base $\C$ of $E$,
then they either project to $q$, or they
wrap around $q$ a non-negative number of times
(i.e. any disk whose boundary is the projected orbit
must intersect $q$ a non-negative number of times).
We have a natural filtration $F_i$ of the chain complex
for $HF_*^{[0,\infty)]}(H_i)$ where the subspace
$F_i$ consists of the subspace of orbits
projecting to $q$ or wrapping around $q$ at most $i$ times.
We then show that $H_*(F_i / F_{i-1}) \cong HF_*(\phi,i)$.
This is because the cylindrical end of $E$
is the same as the cylindrical end for
the fibration $(0,\infty) \times M_\phi$ used to define $HF_*(\phi,i)$
and $F_i / F_{i-1}$
is generated by the orbits which wrap around the base $(0,\infty) \times S^1$
exactly $i$ times around.

Before we prove Theorem \ref{thm:spectralsequence},
we need a preliminary Lemma:
We will first deform our Lefschetz fibration $E$ slightly.
This does not change $SH_*(\widehat{M})$ as the deformation
is compactly supported hence we can use a Moser theorem.
\begin{lemma} \label{lemma:lefschetzfibrationdeformation}
Let $q \in \C$ be a regular value of the Lefschetz fibration
map $\pi : E \rightarrow \C$.
We can deform the Lefschetz fibration $E$ through
a family of Lefschetz fibrations $E_t$,
so that there exists a small neighborhood $U$
around $q$ where the fibration $\pi^{-1}(U)$
is equal to the product fibration
\[\widehat{F} \times U \twoheadrightarrow U.\]
On this region $\theta_F$ splits as a product
$\theta_F + \theta_U$ where
$d\theta_U$ is some volume form on $U$.
This deformation has compact support 
(I.e. the family of $1$-forms $\theta^t_E$
inducing this deformation all agree outside
some compact set).
\end{lemma}
\proof
Because $q$ is a regular value and $\pi$
has only finitely many singularities,
there exists a small neighborhood $V$
where $\pi$ is regular on $\pi^{-1}(V)$.
By using 
parallel transport techniques from \cite[Lemma 8.6, Step 1]{McLean:symhomlef},
we can choose a smooth trivialization
\[\widehat{F} \times V \twoheadrightarrow V\]
of $\pi$ with $\theta_E = \theta_F + \theta_V + b + dR$
where $b$ is some $1$-form which vanishes
when restricted to each fiber $\widehat{F} \times \left\{ x \right\}$.

Let $\rho_t : \C \rightarrow \R$ be a sequence
of functions parameterized by $t \in [0,1]$
such that 
$\rho_t = 1$ near the boundary of $V$ and outside $V$,
$\rho_0 = 1$ everywhere and
$\rho_1$ is equal to $0$ on a smaller neighborhood $U$ of $q$.
We then deform $\theta_E$ through the family
\[ \theta^t_E = \theta_F + \theta_V + \rho_t b + d(\rho_t R)\]
of one forms.
These formulas makes sense because $\rho_t = 1$
outside $V$, so we just define $\theta^t_E$
to be equal to $\theta^t_E$ in this region.

This deformation has compact support
because $b$ and $dR$ are equal to $0$
for $r_F$ large.
Also $\theta^1_E = \theta_F + \theta_V$
which is the product form we want.
This completes our Lemma.
\qed

\proof of Theorem \ref{thm:spectralsequence}.
Let $E,q,U$ be as in Lemma \ref{lemma:lefschetzfibrationdeformation}.
The Lefschetz fibration has one positive end
$C = [1,\infty) \times M_\phi$.
Let $r$ parameterize $[1,\infty)$ and let
$\alpha$ be the contact form on $M_\phi$
such that $\theta_E = r d \vartheta + \alpha$ on this cylindrical end.
Here $\vartheta$ is the pullback via
the fibration map $p_\phi : M_\phi \twoheadrightarrow S^1$
of the angle coordinate of $S^1$.
First of all, we can assume that $U$ is a small disk
centered at $q$
(as we can shrink it). Let $s$ be the radial coordinate
for this disk.
We can also assume that $q$ is disjoint from the cylindrical
end $C$ (either by shifting the cylindrical
end or by moving $q$).
We can also shrink $U$ so that it is also
disjoint from $C$.

Let $H_i$ be a family of Lefschetz Hamiltonians indexed by
$i \in \N$ with the following properties:
\begin{enumerate}
\item \label{item:radialcoordinates}
If $(s,t)$ are radial coordinates for the disk $U$,
$H_i = b_i(s) + H_F$ inside $U$
where $H_F : \widehat{F} \rightarrow \R$ is some admissible Hamiltonian
on the fiber. The $1$-periodic orbits of $H_F$ are non-degenerate.
\item \label{item:conditionnearq}
$b_i'(s) \geq 0$ and is $0$ if and only if $s = 0$.
\item
$b_i'(s)$  tends to $0$ as $i$ tends to infinity.
\item
All the $1$-periodic orbits of $H_i$ are non-degenerate.
For this to be true, $H_i$ really should be an $S^1$
family of Hamiltonians.
\item
$H_i$ is cofinal, which means that
$H_i$ tends to infinity pointwise in the region $\left\{ r > 0 \right\}$
and it tends to $0$ pointwise everywhere else.
\item \label{item:cylindricalcondition}
On the cylindrical end $C$, we have that in the region
$\left\{ 1 < r < 1+1/i \right\}$,
$H_i = h_i(r)$
where $h_i'$ is small so there are no periodic orbits in this region.
We also assume that $h_i' = 1$ in the region $[\frac{1}{i+2},\frac{1}{i}]$.
\end{enumerate}

We also choose a sequence $J_i$
of Lefschetz admissible almost complex structures
such that:
\begin{enumerate}
\item
The pair $(H_i,J_i)$ is regular (i.e. so that 
$HF_*(H_i,J_i)$ is well defined).
For this to be true, we really need an $S^1$ family of almost
complex structures.
\item
Let $j$ be the standard complex structure with
respect to the radial coordinates $(s,t)$ in $U$.
This means that $J_i(\frac{\partial}{\partial s}) = \frac{\partial}{\partial t}$.
We assume that in $\pi^{-1}(U)$,
$J_i = j + J_F$ where $J_F$
is an admissible almost complex structure for $\widehat{F}$.
\item
We require $\pi$ to be $(J,j)$ holomorphic in the
region $1/2i < r < 1/i$.
The reason why we need this condition is to ensure that
we have control over cylinders satisfying the perturbed
Cauchy-Riemann equations passing through this region.
\end{enumerate}

To define symplectic homology,
we don't just need a cofinal family of Hamiltonians,
we also need for two such Hamiltonians $H_j \leq H_i$,
a smooth family of Hamiltonians joining them.
We ensure that each element of this smooth family of Hamiltonians $H^s_{i,j}$
has exactly the same properties listed above as $H_i$ and $H_j$.
The only extra condition is that in $\left\{ 1 < r < 1/i \right\}$
we require $H^s_{i,j} = h^s_i(r)$ where
$\frac{\partial^2}{\partial s \partial r} h^s_i(r) >0.$
We also require similar condition $b^s_i$ inside $U$ where
$b^s_i$ is the function from property (\ref{item:conditionnearq}).
Let $CF_*(H_i,J_i)$ be the chain complex for $HF_*(H_i,J_i)$.
Symplectic homology is the homology of the following complex:
\[\varinjlim_i CF_*(H_i,J_i)\]
where the maps of this directed system $(H_i,J_i)$
come from continuation maps $H^s_{i,j}$ described above.

We can put a filtration on this chain complex as follows:
We first start with filtrations $F^i_k$ on the chain
complex $CF_*(H_i,J_i)$.
The orbits $x$ of $H_i$
which project to $q$ or whose projection to the base winds $0$ times around
$q$ generate a subvector space $F^i_0 \subset CF_*(H_i,J_i)$.

All other projected orbits wrap positively around $q$,
so we define $F^i_k \subset CF_*(H_i,J_i)$
to be the subspace generated by orbits whose
projection winds around $q$ at most $k$ times.
This means that
$\cup_k F^i_k =CF_*(H_i,J_i)$
and $F^i_0 \subset F^i_1 \subset F^i_2 \cdots$. 
We wish to show that if $x$ is an orbit in $F^i_k$,
then $\partial x \in F^i_k$ as well.
This is true if we can show that each solution
$u : \R \times S^1 \rightarrow E$
satisfying the perturbed Cauchy Riemann equations
intersects the fiber $\pi^{-1}(q)$ positively.
This is true because the Hamiltonian vector field
$X_{H_i}$ is tangent to $\pi^{-1}(q)$
and because $\pi^{-1}(q)$ is a holomorphic submanifold.
This means that if we construct the mapping torus
$M_{\phi^1_{H_i}}$ of the symplectomorphism induced
by the $S^1$ family of Hamiltonians $H_i$,
then solutions of the perturbed Cauchy Riemann equations
correspond to holomorphic sections of the holomorphic fibration
\[\R \times M_{\phi^1_{H_i}} \twoheadrightarrow \R \times S^1.\]
The subset $\pi^{-1}(q)$ becomes a holomorphic submanifold of $\R \times M_{\phi^1_{H_i}}$
(i.e. the mapping torus $\R \times M_{\phi^1_{H_i|_{\pi^{-1}(q)}}}$)
of complex codimension $1$.
Any holomorphic curve must intersect this positively.
In particular all holomorphic sections must intersect this positively, and
holomorphic sections correspond to Floer trajectories.
Finally, if a holomorphic section intersects it positively,
then the corresponding solution $u$ of the perturbed Cauchy-Riemann equations
must intersect $\pi^{-1}(q)$ positively.
This implies that $F^i_0 \subset F^i_1 \subset \cdots$
is a filtration.

Solutions of the continuation map equations also must
intersect the fiber $\pi^{-1}(q)$ positively.
Hence the natural continuation maps
$CF_*(H_i,J_i) \rightarrow CF_*(H_j,J_j)$
respect the filtration structure.
This means that it sends $F^i_k$ to $F^j_k$.
Hence on the chain complex
$\varinjlim_i CF_*(H_i,J_i)$ has an induced
filtration $F_0 \subset F_1 \subset \cdots$
where $F_k$ is the direct limit $\varinjlim_i F^i_k$.

The spectral sequence we want is the one induced by
this filtration.
Hence in order to prove our result, we need to show
that the homology of the chain complex
$F_k / F_{k-1}$ is equal to
$HF_{* + 2k -1}(\phi,k)$.
The chain complex $F_k / F_{k-1}$
is the same as the chain complex generated by
orbits wrapping $k$ times around $\pi^{-1}(q)$
and where the differential counts Floer trajectories that do not intersect
this fiber.
We have the Lefschetz fibration
$W_\phi = (0,\infty) \times M_\phi$.
The region $[1,\infty) \times M_\phi \subset W_\phi$
is exactly the same as the cylindrical end $C$ of $E$.
So we create another cofinal family of Lefschetz admissible
Hamiltonians $H_i'$ for $HF_*(\phi,k)$
where $H_i' = H_i$ in the region $r > 1+\frac{1}{2i}$,
and $H_i' = g(r)$ for $r \leq 1+\frac{1}{2i}$
where $g'(r)$ has a very small positive derivative,
so that $H_i'$ has no periodic orbits in this region.
Hence there is a $1-1$ correspondence
between orbits of $H_i'$ that wrap around the base
$k$ times and orbits of $H_i$ that wrap around
$\pi^{-1}(q)$ $k$ times when projected to the base $\C$.
We also choose a Lefschetz admissible almost
complex structure $J_i'$ such that $J_i' = J_i$
in the region $r > 1+\frac{1}{2i}$.
The chain complex of the quotient
$F_k / F_{k-1}$ is generated by these orbits
that wrap around $k$ times.
The differential counts cylinders satisfying
the perturbed Cauchy-Riemann equations joining these orbits.
All these orbits are in the region $r \geq 1+1/i$,
and also Lemma \ref{lemma:lefschetzminimumprinciple}
ensures that any cylinder satisfying the perturbed
Cauchy-Riemann equations stays in this region as well.
The reason why Lemma \ref{lemma:lefschetzminimumprinciple} works here is because
these cylinders do not intersect the fiber $\pi^{-1}(q)$,
hence these are subsets of the fibration
\[E \setminus \pi^{-1}(q) \twoheadrightarrow \C \setminus q = C_*.\]

The same argument also ensures that
all the orbits and cylinders satisfying the perturbed Cauchy-Riemann
equations in $W_\phi$ stay inside the region $r \geq 1/i$ as well.
We then have $F^i_k / F^i_{k-1}$
is chain isomorphic to
$CF^{\alpha_k}_{* + 2k -1}(H_i',J_i')$,
hence $H_*(F^i_k/F^i_{k-1}) = HF^{\alpha_k}_{* + 2k -1}(H_i',J_i')$.
The shift in grading comes from the fact that we
have two different trivializations of the tangent bundle of
the base $\C_*$. The first trivialization
comes from embedding $\C_*$ in $\C$ and choosing
the standard trivialization of $T\C$.
The other trivialization comes from identifying
$\C_*$ with $\C / \Z$ where the $\Z$ action
is generated by the map $(x + iy) \rightarrow x + i(y+1)$
and then trivializing $T\C$ in the standard way.
Any such trivialization of $T\C_*$
combined with an appropriate sequence of
trivializations of the canonical bundle of
the fiber $\widehat{F}$
induces a trivialization of the canonical bundle of our respective
Lefschetz fibrations.
This is done by using the
splitting of the tangent bundle into
horizontal and vertical subspaces
(See \ref{section:equivariantfloer} earlier).

The continuation maps between these Floer homology groups
are the same
as well using a similar arguments,
hence
\[H_*(F_k / F_{k-1}) = \varinjlim_i H_*(F^i_k/F^i_{k-1}) =\]
\[\varinjlim_i HF^{\alpha_k}_{*+2k-1}(H_i',J_i')
= HF_{*+2k-1}(\phi,k).\]

This gives us our spectral sequence.
\qed

\subsection{Group actions} \label{section:groupactions}

In this section we will prove
Lemma \ref{thm:groupactions} which says:
Suppose that the coefficient field $\K$
has characteristic $0$ or characteristic $p$
where $p$ does not divide $k$,
then there exists a $\Z / k\Z$ action
$\Gamma$ on $HF_*(\phi^k,1)$
such that
$HF_*(\phi,k) \cong HF_*(\phi^k,1)^\Gamma$.

\proof
Let $\pi_\phi : W_\phi \twoheadrightarrow \R \times S^1$
be equal to the Lefschetz fibration $(0,\infty) \times M_\phi$.
Then there is a natural $k$ fold covering map:
$p_k : W_{\phi^k} \twoheadrightarrow W_\phi$.
Here $W_{\phi^k}$ is obtained as the pullback
bundle of $W_\phi$ via the map
$c : (0,\infty) \times S^1 \twoheadrightarrow (0,\infty) \times S^1$
where $c(s,t) = (s,kt)$ where
we view $S^1$ as $\R / \Z$.

Let $\beta_k \subset H_1(W_\phi)$ be the set of
$H_1$ classes represented by loops which project
down to loops in $(0,\infty) \times S^1$ which wrap
around $S^1$ $k$ times.
Then $({p_k})_*^{-1}(\beta_k)$ is the generated by loops
in $W_{\phi^k}$ which wrap around the base once.
If we have a Lefschetz admissible pair
$(H_i,J_i)$ on $W_\phi$, then the preimage
$(\widetilde{H_i},\widetilde{J_i})$ of $(H_i,J_i)$ under the covering map $p_k$
is also Lefschetz admissible.
Also if $(H_i,J_i)$ is regular for orbits wrapping $k$ times around, then
so is $(\widetilde{H_i},\widetilde{J_i})$.
Here is a very brief sketch of why this statement is true:
Non-degenerate
orbits of $H_i$ wrapping $k$ times around the base
lift to non-degenerate orbits wrapping once around the base
(there is a choice of $k$ lifts).
Also if $x,y$ are such orbits of $H_i$ and
$\widetilde{x},\widetilde{y}$ are choices of lifts of $x$ and $y$,
then the moduli space of cylinders joining
$\widetilde{x}$ and $\widetilde{y}$ is a clopen component
of the moduli space of cylinders joining
$x$ and $y$.
Hence if $J_i$ is regular for this moduli space then $\widetilde{J}$ must also be regular
for this lifted moduli space
because the linearized $\bar{\partial}$ operator is exactly the same.

The deck transformations of the covering map $p_k$
preserve $(\widetilde{H_i},\widetilde{J_i})$.
These transformations induce a $\Z / k\Z$ action on
$CF_*^{p_k^* \beta_k}(\widetilde{H_i},\widetilde{J_i})$.
Let $\Gamma : \Z / k\Z \rightarrow \text{Hom}(CF_*^{p_k^* \beta_k}(\widetilde{H_i},\widetilde{J_i}))$
be the induced action on the chain complex.
The quotient complex:
$CF_*^{p_k^* \beta_k}(\widetilde{H_i},\widetilde{J_i}) / (\Z / k \Z)$
is canonically isomorphic to the chain complex
$CF_*^{\beta_k}(H_i,J_i)$.
The $\Z / k\Z$ action on the set of orbits of $\widetilde{H_i}$
(which generate the vector space $CF_*^{p_k^* \beta_k}(\widetilde{H_i},\widetilde{J_i})$)
is free, hence we can use
work from \cite[Proposition 3G.1]{Hatcher:algebraictopology}
to construct a transfer map from
$HF_*^{p_k^* \beta_k}(\widetilde{H_i},\widetilde{J_i})$
to $HF_*^{\beta_k}(H_i,J_i)^\Gamma$
(the subgroup consisting of elements invariant
under the action of $\Gamma$).
If the coefficient field $\K$ has characteristic $0$
or characteristic $p$ not dividing $k$, then
this map is an isomorphism.
The continuation maps are compatible with the group actions
as well, hence taking direct limits gives us our result.
This proves our lemma.
\qed

\subsection{A long exact sequence between Floer homology groups}
\label{section:floerlongexactsequence}

We will prove theorem \ref{thm:floerlongexactsequence}.
Here the statement of this theorem:
For any symplectomorphism $\phi$, we have a long exact sequence
\[\rightarrow HF^i(\phi) \rightarrow HF_i(\phi,1) \rightarrow
HF^{i-1}(\phi) \rightarrow. \]
Before we prove this Theorem, we need a preliminary Lemma.
This is a correspondence between maps into a particular Lefschetz
fibration and maps into its fiber.
Suppose we have a non-degenerate symplectomorphism
$\phi : F \rightarrow F$, and let $W_\phi$
be the mapping cylinder constructed as follows:
Let $\widetilde{S} = (0,\infty) \times [0,1]$  with the symplectic form
$ds \wedge dt$ where $s$ parameterizes the first interval
and $t$ parameterizes the second one.
We take $S := \widetilde{S} / \sim$ where $\sim$ identifies $(t,0)$
with $(s,1)$. The symplectic form $ds \wedge dt$ descends to $S$.
We define $\widetilde{W} := \widetilde{S} \times F$
with the natural product symplectic structure
and we define $W := \widetilde{W} / \sim$
where $\sim$ identifies $(t,0,f)$ with $(t,1,\phi(f))$
where the first two coordinates are the natural coordinates parameterizing
$\widetilde{S}$ and the third one is a point $f \in F$.
Again the symplectic form descends.
Let $p_1,p_2$ be the projection maps from $\widetilde{W}$
to $\widetilde{S}$ and $F$ respectively.
The map $p_1$ descends to a map $\pi_W$ from $W$ to $S$.

We have a natural connection on $\pi_W$ given by the $\omega$-orthogonal
plane distribution to the fibers of $\pi_W$.
If $H_t$ is any Hamiltonian on $S$ then the Hamiltonian
vector field of $\pi_W^*H_t$ at a point $x$ is the unique horizontal lift
at $x$ of $X_{H_t}$ in $TS_{\pi_W(x)}$.
Let $H_t$ be an $S^1$ family of Hamiltonians on $S$
such that it has a non-degenerate orbit $l : S^1 \rightarrow S$
where $l$ is the natural injection $S^1 = \{1\} \times S^1 \hookrightarrow S$.
We will also suppose that $\phi$ is a non-degenerate symplectomorphism.
Let $J_T$, $T \in [0,1]$ be a family of almost complex structures parameterized
by $[0,1]$ on $F$ such that $\phi_* J_0 = J_1$.
Let $j_S$ be a complex structure on $S$.
From these we can construct a new almost complex structure on
$T\widetilde{W} = T\widetilde{S} \oplus TF$ given by $j_S \oplus J_t$
where $t$ is the coordinate parameterizing $[0,1]$ in $\widetilde{S} = (0,\infty) \times [0,1]$.
Conversely suppose we have some almost complex structure $J$
on $W$ which is compatible with the symplectic form
and makes $\pi_W$ $(J,j_S)$ holomorphic then
we can construct a family of almost complex structures
$J_t$ on $F$ such that $\phi_* J_0 = J_1$.
The point is that $J$ must split up as $J_S \oplus J_t$
after we pull it back to $T\widetilde{W} = T\widetilde{S} \oplus TF$
because the projection map to $\widetilde{S}$ is $(J,j_S)$ holomorphic
and compatible with the symplectic form.

\begin{lemma} \label{lemma:modulispacebijection}
There is a natural 1-1 correspondence $\Psi_p$ between $1$-periodic orbits
of $\pi_W^* H_t$ which project to $l$ and fixed points of $\phi$.
Let $\D$ be a small disk around $l(0)$.
If we have two such $1$-periodic orbits $x,y$ of $\pi_W^* H_t$
then there is a natural 1-1 correspondence between
smooth maps $u : \R \times S^1 \rightarrow W$
connecting $x,y$
such that
\begin{enumerate}
\item $\pi_W(u(0,0)) \in \D$.
\item The loop $\pi_W(u(0,t))$ is homotopic to $l$
\end{enumerate}
and pairs of smooth maps $u_1 : \R \times [0,1] \rightarrow F$,
$u_2 : \R \times S^1 \rightarrow S$
such that
\begin{enumerate}
\item $u_1(s,1) = \phi(u_1(s,0))$.
\item $u_1(s,t)$ converges to $\Psi_p(x)$ as $s$ goes to $-\infty$.
and converges to $\Psi_p(y)$ as $s$ goes to $+\infty$.
\item $u_2(s,t)$ converges to $l(t)$ as $s$ goes to $\pm \infty$.
\item $u_2(0,0) \in \D$.
\item The loop $u_2(0,t)$ is homotopic to $l$.
\end{enumerate}
Also the Cauchy-Riemann operator
$\partial_{J,H}(u) := (Du + \pi_W^*X_{H_t} \otimes ds)^{(0,1)}$
maps naturally to the sum
$\partial_J(u_1) + \partial_{J,H}(u_2)$.
\end{lemma}
\proof of Lemma \ref{lemma:modulispacebijection}.
Let $q : S^1 \rightarrow W$ be a $1$-periodic orbit
of $H_t$ which maps to $l$.
Let $Q_W : \widetilde{W} \twoheadrightarrow W$,
$Q_S : \widetilde{S} \twoheadrightarrow S$ be the natural quotient maps.
We can lift $q$ to a map $q : [0,1] \rightarrow \widetilde{W}$
which is the Hamiltonian flow of $Q_W^* \pi_W^* H_t = (p_1)^* Q_S^* H_t$.
We define $\Psi_p(q)$ to be the point $p_2(q(0)) \in F$.
We have a natural inverse as follows:
If $x$ is a fixed point then we take the line given by the inclusion
\[v : [0,1] = \{1\} \times [0,1] \times \{x\} \hookrightarrow \widetilde{W}
 = (0,\infty) \times [0,1] \times F.\]
This is a flow line for $(p_1)^* Q_S^* H_t$
and it satisfies $\phi(v(0)) = v(1)$ hence it descends to a $1$-periodic orbit
of $H_t$.
This is exactly the inverse of $\Psi_p$.

We now wish to find a natural $1$-$1$ correspondence between
maps $u$ and pairs $(u_1,u_2)$.
We define $u_2$ as $\pi_W \circ u$ this satisfies $u_2(0,0) \in \D$.
Let $\widetilde{W}'$ be the universal cover of $W$.
This is symplectomorphic to the product $(0,\infty) \times \R \times F$
with the product symplectic form $ds \wedge dt + \omega_F$.
The quotient map $\pi'_W : \widetilde{W}' \twoheadrightarrow W$
is induced by the quotient of the $\Z$ group action where $1 \in \Z$
is the map:
$(s,t,f) \rightarrow (s,t+1,\phi(f))$.
Choose some lift $\widetilde{\D}$ of $\D$ in $\widetilde{W}'$.
We now replace $u$ with $\widetilde{u} : \R \times [0,1] \rightarrow W$
which is the composition $u \circ p_{S'}$ where
$p_{S'}$ is the natural projection map $\R \times [0,1]$ to
$\R \times S^1$ identifying $(t,0)$ with $(t,1)$.
Because the domain of $\widetilde{u}$ is contractible, we have
a unique lift of $\widetilde{u}$ to $\widetilde{u}' : \R \times [0,1] \rightarrow \widetilde{W}'$
such that $\widetilde{u}'(0,0) \in \widetilde{\D}$.
This is because the disk $\widetilde{\D}$ is small enough
so that the $\Z$ action on $\widetilde{W}'$ sends $\widetilde{\D}$
to other sets disjoint from $\widetilde{\D}$.
Note that the projection maps $p_1$ an $p_2$ extend to $\widetilde{W}'$
if we view $\widetilde{W}$ as some subset of $\widetilde{W}'$.
We define $u_1$ to be $p_2 \circ \widetilde{u}'$.
Because the loop $\pi_W(u(s,t))$ is homotopic to $l$ for each $s$,
we get that $u_1(s,0) = \phi(u_1(s,1))$.

Suppose we are now given $u_1$ and $u_2$, we wish to reconstruct
$u$ from these maps.
Let $\widetilde{u_2}$ be the composition $u_2 \circ p_{S'}$.
Choose a unique lift $\widetilde{u_2}'$ of $\widetilde{u_2}$ into the universal
cover $\R \times \R$ such that $\widetilde{u_2}'(0,0) \in p_1(\widetilde{\D})$.
We have a natural map from $\R \times [0,1]$ into $\widetilde{W}'$
given by $(u_1,\widetilde{u_2}')$.
This projects to the map $u$.

We have $\partial_{J,H}(u) = (Du + X_{\pi_W^* H_t} \otimes ds)^{(0,1)}$,
$\partial_J(u_1) = (Du_1)^{(0,1)}$
and $\partial_{J,H}(u_2) = (Du_2 + X_{H_t} \otimes ds)^{(0,1)}$.
Also for any vector $X$, $(\pi_W)_*(X)$ is equal to
$(\pi_W)_*(X^h)$ where $X^h$ is the horizontal component of $X$.
We have
\[(\pi_W)_*(Du + X_{\pi_W^* H_t} \otimes ds)^{(0,1)}(Y) =\]
\[(\pi_W)_*\left( Du(Y)^h + X_{\pi_W^* H_t}^h(ds(Y)) +(J \circ Du(jY))^h +
 JX_{\pi_W^* H_t}^h(ds(jY)) \right).\]
We have $X_{\pi_W^* H_t} = \widetilde{X}_{H_t}$
where $\widetilde{X}_{H_t}$ at $x$ is the unique horizontal lift
of $X_{H_t}$ at $\pi_W(x)$.
Also $\pi_W$ is $(J,j_S)$-holomorphic.
Hence we have \[\pi_W Du(Y)^h = Du_2(Y),
\pi_W(X_{\pi_W^* H_t}^h) = X_{H_t},\]
\[(J \circ Du(jY)(Y))^h = j_S \circ Du_2(jY).\]
$ \pi_W(JX_{\pi_W^* H_t}^h(ds(jY))) = j_S \circ X_{H_t} (ds(jY))$.
This implies
$(\pi_W)_* \partial_{J,H}(u) = \partial_{J,H} u_2$.

The Cauchy-Riemann operator for the map $\widetilde{u}$
is exactly the same as the one for $u$
(you just pull back everything via the map $p_{S'}$).
Also if we consider the lift $\widetilde{u}'$
of $\widetilde{u}$ as described earlier,
we still have exactly the same operator.
This is because $\widetilde{W}'$ is a cover
of $W$ such that its projection map is holomorphic.
The almost complex structure on $W$ is $(j_S,J_t)$
on the tangent space $\R \times \R \times TF$
of $\widetilde{W}'$ at the point $(s,t,f)$.
We defined $T^v\widetilde{W}'$ as the tangent
spaces $TF$ in $\R \times \R \times TF$.
These are the vertical tangent spaces of the fibration $p_1$.
We have that: 
\[(p_2)_*(D{\widetilde{u}'} + X_{p_1^* H_t} \otimes ds)^{(0,1)}(Y)
= (p_2)_*(D{\widetilde{u}'} + X_{p_1^* H_t} \otimes ds)^{(0,1)}(Y)^v.\]
Because $X_{p_1^* H_t}$ vanishes on $T^v\widetilde{W}'$,
we get that this is equal to:
\[(p_2)_* ((D_{\widetilde{u}'}(Y))^v)^{(0,1)} =
((D{\widetilde{u}'}(Y))^v + J_t \circ (D_{\widetilde{u}'}(jY))^v  ).\]
Because $(p_2)_* D_{\widetilde{u}}^v(Y) = D(u_1)(Y)$, we get:
\[(p_2)_* \partial_{J,p_1^*H}(\widetilde{u}')(Y) =
(D(u_1)(Y)) + J_t \circ D(u_1)(jY) = (\partial_J(u_1))(Y).\]

Hence our correspondence between $u$ and $(u_1,u_2)$ also gives
us a natural map between Cauchy-Riemann operators on $u$
and $(u_1,u_2)$ respectively.
\qed

Because we have a nice correspondence between between maps $u$ and
$(u_1,u_2)$ and their respective Cauchy-Riemann operators,
we have that the moduli space of maps $u$ satisfying
the perturbed Cauchy-Riemann equations is the same
as the moduli space of $J_t$ holomorphic maps $u_1$
and maps $u_2$ satisfying the perturbed Cauchy-Riemann equations.
If we have regularity for one moduli space then we have regularity
for the other and we can ensure that their orientations coincide.
Because the maps $u_2$ satisfying the perturbed Cauchy-Riemann
equations joining $l$ with $l$ have energy zero,
they must map to the constant loop $l$.
Hence we actually have a natural correspondence
between maps $u$ satisfying the perturbed Cauchy-Riemann equations
and $J_t$-holomorphic maps $u_1$.

Let $J_T,T \in [0,1]$ be a smooth family of almost complex structures
compatible with the symplectic form on $F$
such that $\phi_* J_0 = J_1$.

\proof of Theorem \ref{thm:floerlongexactsequence}.
Let $\phi'$ be a standard perturbation for
$\phi$ as in definition \ref{label:standardperturbation}.
Let $r_F$ be the cylindrical coordinate for $\widehat{F}$.
Using this symplectomorphism $\phi'$ we will first carefully
construct a Lefschetz fibration with one positive
and one negative end whose monodromy
is almost equal to $\phi'$.

The reason why it cannot be exactly equal to $\phi'$
is that the parallel transport maps are equal
to the identity map for $r_F$ large, but this isn't
true for $\phi'$.
Instead we produce a new symplectomorphism $\phi''$ as follows:
There exists a small constant $\epsilon > 0$ and
a large constant $R > 0$ such that in the region
$r_F \geq R$, $\phi'$ is equal to the symplectomorphism
$\phi^1_{\epsilon r_F}$.
Choose a smooth function $l : [R,\infty) \rightarrow \R$
with the following properties:
\begin{enumerate}
\item $l(r_F) = \epsilon r_F$ for $r_F$ near $R$
\item $0 < l'(r_F) \leq \epsilon$.
\item $l'(r_F) = 0$ for $r_F \geq R+1$
\end{enumerate}
Let $\phi''$ be the new symplectomorphism
equal to $\phi'$ in the region $r_F < R$
and equal to $\phi^1_{l(r_F)}$ in the region $r_F \geq R$.
In section \ref{section:lefschetzfibrations},
we constructed a mapping torus $M_\phi$ with
an explicit contact form equal to
$Cdt + \theta_F + dG$.
The good thing about this contact form is that
the monodromy map around the  base $S^1$
is exactly the same as $\phi$.
Using this construction, let $M_{\phi''}$
be the mapping torus of $\phi''$ with contact
form $\alpha_{\phi''}$ whose monodromy map is
exactly equal to $\phi''$.

Using this contact form $\alpha_{\phi''}$
we can construct a Lefschetz fibration
$W_{\phi''} = (0,\infty) \times M_{\phi''}$
with associated one form
$\theta_{\phi''} = rd\vartheta + \alpha_{\phi''}$.

We will now construct a Lefschetz admissible Hamiltonian
$H$ of slope $\beta$ as follows:
Let $L : W_{\phi''} \rightarrow \R$
be equal to $0$ in the region $\left\{ r_F \leq R \right\}$
and equal to $\epsilon r_F - l(r_F)$ in the region
$\left\{ r_F \geq R \right\}$.
Let $\kappa : (0,\infty) \rightarrow \R$
satisfy:
\begin{enumerate}
\item $\kappa'' \geq 0$ and $\kappa' > 0$.
\item $\kappa'(r) = \beta$ for $r$ large.
\item $\kappa'(r)$ is small for $r$ small.
\item $\kappa'(1) = 2\pi$, $\kappa''(1)>0$.
\end{enumerate}
Let $K : (0,\infty) \times S^1 \rightarrow \R$
be equal to $\kappa(r)$.
Let $H = \pi_\phi^* K + L$.
Any $1$-periodic orbit of $H$ with $H_1$ class in $\beta_1$
projects isomorphically to the circle
$\left\{ r = 1 \right\}$.
Because $H$ is not time dependent, we have $S^1$ families of orbits.
Choose a Morse function $f$ on $S^1 \cong \left\{ r=1 \right\}$ with one maximum and one minimum.
We can use work from \cite{FloerHofer:SymhomI},
combined with the Morse function on each
orbit to perturb $H$ to a non-degenerate Hamiltonian.
In fact, we can perturb $K$ to a time dependent
function $K_t$ so that
$H_t := \pi_\phi^*K_t + L$
has non-degenerate orbits because the symplectomorphism $\phi'$
is non-degenerate.
The Hamiltonian $H_t$ has the following property:
the $1$-periodic orbits of $H_t$ are exactly $1$-periodic orbits
of $H$ that start and end at the critical points of the function $\pi_{\phi''}^*f$.
Let $J$ be an almost complex structure
which is Lefschetz admissible and invariant
under translations in the $r$ direction.
We also assume that $\pi_{\phi''}$
is $(J,j)$ holomorphic.


We view the base $(0,\infty) \times S^1$ as a symplectic manifold
with symplectic form $dr \wedge dt$.
The Hamiltonian $K_t$ has two orbits that wrap around the $S^1$ factor of
$\R \times S^1$ once. One orbit starts and ends at the maximum of the Morse
function $f$ and the other starts and ends at the minimum.
Let $o_m,o_M : S^1 \rightarrow (0,\infty) \times S^1$ be these two orbits
corresponding
to the minimum and maximum of the Morse function $f$
respectively.
Let $p_m,p_M$ be the starting points of the orbits $o_m$ and $o_M$
respectively (these are fixed points of the time $1$ Hamiltonian symplectomorphism
$\phi^1_{H_t}$).
Every fixed point of $\phi^1_{H_t}$ whose associated orbit
wraps around $(0,\infty) \times S^1$ once maps to either
$p_m$ or $p_M$ via $\pi_{\phi''}$.
Define $F_m := \pi_\phi^{-1}(p_m)$ and $F_M := \pi_\phi^{-1}(p_M) \cong \widehat{F}$.
The time $1$ flow of $H_t$ sends $F_m$ to $F_m$ and $F_M$ to $F_M$
and both these maps are exactly equal to $\phi'$
after identifying $F_m$ and $F_M$ with $\widehat{F}$.
This means that there is a bijection between
fixed points of $\phi'$ and fixed points of $\phi^1_{H_t}$
that project to $p_m$. Similarly there is a bijection
between fixed points of $\phi'$ and fixed points
of $\phi^1_{H_t}$ that project to $p_M$.
If $x$ is a fixed point of $\phi'$,
we write $x_m,x_M$ for the corresponding fixed points
of $\phi^1_{H_t}$ that project to $p_m$ and $p_M$ respectively.
We also write $x_m^o,x_M^o$ for the corresponding orbits.

Let $x^1,x^2$ be two fixed points of $\phi$
such that they project down to the fixed
point $p_m$.
Using the orbit $o_m : S^1 \rightarrow (0,\infty) \times S^1$
we can construct a family of almost complex structures $J_t$
on $\widehat{F} \cong F_m$ from $J$ as in
the statement before Lemma \ref{lemma:modulispacebijection}.

Let $\M$ be the moduli space of
$J_t$ holomorphic maps 
$u : \R \times [0,1] \rightarrow \widehat{F}$
where $\phi'(u(1,t)) = u(0,t)$
joining $x^1$ and $x^2$.
Let $\M_m$ be the moduli space of
maps $u : \R \times S^1 \rightarrow \R \times M_\phi$
satisfying the perturbed Cauchy Riemann equations
for $H$ joining $x^1_m$ and $x^2_m$.
We have a natural bijections
\[\M \rightarrow \M_m\]
from \ref{lemma:modulispacebijection}.
Also we can ensure that both are regular and
that the orientation of these moduli spaces are the same.
Similar reasoning ensures that we also have a bijection
between $\M_M$ and $\M$.

Let $\partial$ be the differential for $HF^*(\phi',1)$.
Let $x,y$ be fixed points of $\phi'$.
Regularity of $(K_t,j)$ ensures that
there is no cylinder satisfying Floer's
equations for $H_t$ starting at $x_M$ and ending at $y_m$.
This is because any cylinder $u$ satisfying such equations
projects to another cylinder $u'$ satisfying Floer's
equations for $K_t$ on $(0,\infty) \times S^1$.
But the cylinder $u'$ starts at an orbit of
index $0$ (the orbit $o_M$) and ends at an orbit of index $1$
(the orbit $o_m$). This is impossible by regularity.
Hence we have an increasing filtration $F_M \subset F_M \bigoplus F_m = C_*(H_t)$
where $F_M$ consists of fixed points of the form $x_M$
and $F_m$ consists of fixed points of the form $x_m$.
This means that the differential is of the form:
\[
\partial =
\left( \begin{array}{cc}
\partial_{F_m} & 0 \\
\partial_{m,M} & \partial_{F_M} \end{array} \right).\]

Let $\text{Fix}(\phi')$ be the set of fixed points of
$\phi'$. 
There is a bijection between $\text{Fix}(\phi')$
and orbits of the form $x_m$ given by the map $x \rightarrow x_m$.
Similarly there  is a bijection between $\text{Fix}(\phi')$ 
and orbits of the form $x_M$.
Using these bijections and the fact that we proved
$\M = \M_m = \M_M$ for all orbits $x^1,x^2$,
we have that the chain complexes $(F_M,\partial_{F_M})$
and $(F_m,\partial_{F_m})$ are chain isomorphic to
the chain complex for $HF_*(\phi) = HF_*(\phi')$.
Using this fact and this filtration we get our
long exact sequence.
\qed

\section{Applications of our spectral sequence} \label{section:applications}


We will first prove Corollary \ref{corollary:nonvanishingfloerhomology}.
Here is a statement of this corollary:
Suppose that $\phi : F \rightarrow F$ is a symplectomorphism
such that it is
obtained by one or more stabilizations to the identity map
$\text{id} : F' \rightarrow F'$. Then if 
the Euler characteristic is odd,
then $HF^*(\phi^k,\Q) \neq 0$ for infinitely
many $k$.
From now on (from the comment at the end of Section \ref{section:floerhomologydefn})
we will assume that $\phi$ is a compactly supported
symplectomorphism $\phi : \widehat{F} \rightarrow \widehat{F}$.
We first some preliminary Lemmas.

\begin{lemma} \label{lemma:floerhomologynontrivial}
Suppose that $HF_*(\phi,k,\Q) \neq 0$, then
$HF^*(\phi^k,\Q)$ is non-trivial.
\end{lemma}
\proof of Lemma \ref{lemma:floerhomologynontrivial}.
Theorem \ref{thm:groupactions} tells us that there is a $\Z / k\Z$
action on $HF_*(\phi^k,1,\Q)$ whose fixed points
give us $HF_*(\phi,k,\Q)$.
This means that if $HF_*(\phi,k,\Q)$
is non-zero, then so is $HF_*(\phi^k,1,\Q)$.
Then the long exact sequence from Theorem \ref{thm:floerlongexactsequence}
tells us that $HF^*(\phi^k,\Q)$ is non-zero.
\qed

\begin{lemma} \label{lemma:oddrank}
Let $E^*_{*,*}$ be a spectral sequence (with coefficients
in some field $\K$).
Suppose for some $r \geq 0$,
the total rank of $E^r_{*,*}$ is odd (resp. even),
then the total rank of $E^\infty_{*,*}$ is odd (resp. even).
\end{lemma}
\proof of Lemma \ref{lemma:oddrank}.
This is done by induction.
Suppose for some $R \geq r$, $E^R_{*,*}$ has odd rank.
Then there is a differential $\partial$ on $E^R_{*,*}$
such that $E^{R+1}_{*,*} = H_*(E^R_{*,*},\partial)$.
This must have odd rank, because the homology
of an odd dimensional chain complex is odd dimensional.
The reason why this is true is because:
if the rank of the image of $\partial$ is odd,
then the rank of the kernel of $\partial$ is 
even by the first isomorphism theorem.
This means that the rank of $H_*(E^R_{*,*},\partial)$ is an odd number
minus an even number which is odd.
Similar reasoning ensures that if the rank of the image
of $\partial$ is even, then the homology has odd rank.

Exactly the same reasoning ensures that if $E^r_{*,*}$
has even rank, then so does $E^R_{*,*}$.
This proves the Lemma.
\qed

\begin{lemma} \label{lemma:evendimensionalhomology}
The rank of $HF_*(\phi,k,\Q)$ is even.
\end{lemma}
\proof of Lemma \ref{lemma:evendimensionalhomology}.
The idea here is to construct an explicit chain
complex of even rank for $HF_*(\phi,k,\Q)$.
Then the same reasoning as in the previous
Lemma \ref{lemma:oddrank} gives us our result.
Let $W_\phi := \R \times M_\phi$ be the Lefschetz
fibration associated to $\phi$.
Let $H : W_\phi \rightarrow \R$ be a Lefschetz admissible
Hamiltonian which is {\it time independent}.
Using similar methods to \cite[Theorem 3.1]{HoferSalamon:FloerNovikov},
a generic such $H$ has non-degenerate orbits.
This means that for all fixed points $x$ of $\phi^1_H$,
the linearized return map 
\[D_x \phi^1_H : T_x W_\phi \rightarrow T_x W_\phi\]
has at most one eigenvalue equal to $1$.
Let $\beta_k \subset H_1(W_\phi)$ be equal
to ${\pi_\phi}_*^{-1}(l_k)$ where $l_k$
is the homology class represented by a loop
wrapping $k$ times positively around the $S^1$ factor of the base
$\R \times S^1$.
All of the $1$-periodic orbits whose $H_1$
class is in $\beta_k$ are non-trivial.
Choose a Morse function $f$ on each $1$-periodic orbit $o$
of $\phi^1_H$ whose $H_1$ class is
in $\beta_k$ with one maximum and one minimum.
We can do this because each of these orbits are non-trivial.
We can perturb $H$ to a time dependent Hamiltonian
\[H' : S^1 \times W_\phi \rightarrow \R\]
using these Morse functions so that
for each $1$-periodic orbit $x$ of $H$ in the class
$\beta_k$, we get two non-degenerate periodic orbits $x_m,x_M$
of $H'$ in the class $\beta_k$ corresponding
to the minimum $m$ and maximum $M$ of the Morse function $f$
(see \cite{CieliebakFloerHoferWysocki:SymhomIIApplications}).
This implies that the chain complex
$CF_*^{\beta_k}(H',J)$ is even dimensional, which
means that $HF_*^{\beta_k}(H',J)$ is even dimensional.
We can construct a cofinal family $H_i'$ of such Hamiltonians
so that $HF_*(\phi,k,\K) = \varinjlim_i HF_*^{\beta_k}(H_i',J)$.
By Lemma \ref{lemma:alternativefloerdefinition},
there exists an $i$ such that
$HF_*(\phi,k,\K) = HF_*^{\beta_k}(H_i',J,\K)$.
This implies that the homology group $HF_*(\phi,k,\K)$
has even rank.
\qed

\proof of Corollary \ref{corollary:nonvanishingfloerhomology}.
In view of Lemma \ref{lemma:floerhomologynontrivial}
we have to prove that $HF_*(\phi,k,\Q)$ is non-trivial
for infinitely many $k$.
We have a Lefschetz fibration:
\[\pi_{\text{prod}} : \widehat{F'} \times \C \twoheadrightarrow \C\]
which is the natural projection whose monodromy map
is the identity $\text{id} : \widehat{F'} \rightarrow \widehat{F'}$.
By \cite[Proposition 2]{Oancea:kunneth}, we get that
$SH_*(\widehat{F'}) = 0$.
Positive stabilization
is some operation on a Lefschetz fibration
which does not change the symplectomorphism
type of the total space
(see Theorem \ref{thm:stabilizationtheorem} in the appendix).
In particular, any sequence
of positive stabilizations of $\pi_{\text{prod}}$
gives us a new Lefschetz fibration:
\[\pi_E : E \twoheadrightarrow \C\]
such that $E$ is symplectomorphic
to $\widehat{F'} \times \C$.
Hence $SH_*(E) = 0$.
We will assume that the monodromy
map of $E$ is equal to $\phi$.

Suppose for a contradiction
that there exists a $K>0$ such that
that $HF_*(\phi,k,\Q) = 0$ for all $k \geq K$.
Then by Theorem \ref{thm:spectralsequence},
there exists a spectral sequence
converging to $SH_*(E) = 0$
with $E^1$ page:
\[E^1_{0,q} = H^{n-*}(E)\]
And for $p > 1$,
\[E^1_{p,q} = HF_{q-p+1}(\phi,q).\]
For $p<0$, 
\[E^1_{p,q} = 0.\] 
The total rank
$\bigoplus_{p,q} E^1_{p,q}$ is finite and odd.
This is because the rank of $E^1_{0,q} = H^{n-*}(E,\Q)$
is odd and the rank of $E^1_{p,q}$ is even
for $p \neq 0$ and $E^1_{p,q} = 0$ for $|p|,|q| \gg 0$.
The reason why the rank of $H^{n-*}(E,\Q)$
is odd is because the Euler characteristic of $E$ is odd.
By Lemma \ref{lemma:oddrank}, this implies that
that the rank of $\bigoplus_{p,q} E^\infty_{p,q} \cong SH_*(E)$
has odd rank.
But $SH_*(E)$ has even rank (equal to $0$).
\qed

The above proof actually tells us the following fact:
For any Lefschetz fibration $E$, if
the rank of $SH_*(E) \text{ mod } 2$ is different
from the rank of $H^*(E) \text{ mod } 2$,
then $HF^*(\phi^k)$ does not vanish for infinitely many $k$
where $\phi$ is the monodromy map.
Similar methods also show that if the rank of $SH_*(E)$
is infinite, then $HF^*(\phi^k)$ does not vanish for infinitely
many $k$. For instance $SH_*(T^*M)$ has infinite rank
if $M$ is simply connected, hence the monodromy
map of any Lefschetz fibration symplectomorphic to $T^*M$
has the above property.


\section{Appendix A: Stabilization}

Let $M$ be an exact symplectic manifold with boundary.
Let $d\lambda$ be the associated symplectic form on $M$.
A {\it Weinstein function} $f : M \rightarrow \R$
is a function such that $-i(X_f)\lambda > 0$ away from the critical points of $f$.
A {\it Weinstein cobordism} is an exact symplectic manifold $M$
as above with a Weinstein function $f : M \rightarrow [c,d]$
where $\partial M = \partial_- M \sqcup \partial_+ M$ with
$\partial_- = f^{-1}(c)$ and $\partial_+ = f^{-1}(d)$ and such that
$c$ and $d$ are regular values of $f$.
We will assume that $M$ is compact.

\subsection{Attaching a Weinstein $n$-handle}

We will first describe how adding Dehn twists
corresponds to attaching $n$-handles.
We will do this in two parts.
We will first describe the $n$-handle carefully
and attach it to the Lefschetz fibration
creating a Liouville domain.
We will then deform this Liouville domain
so that the Lefschetz fibration structure
extends over the handle in such a way that
we have a new critical point of our Lefschetz
fibration whose vanishing cycle corresponds
to the attaching Legendrian sphere.

\begin{defn} \label{defn:compactlefschetzfibration}
Let $\pi : E \twoheadrightarrow S$ be a Lefschetz fibration as in Definition
\ref{defn:lefschetzfibration}.
As in the definition, there is a subset $E_h \subset E$
which is symplectomorphic to a product $S \times \partial F \times [1,\infty)$.
Let $E^0_h \subset E^h$ be the interior of this set
(I.e. the subset $E_h \setminus \left(  S \times \partial F \times \left\{ 1 \right\}\right)$).
Let $A$ be the union of all the positive and negative ends
of $E$. Let $A^0$ be equal to the interior of $A$.
Let $\bar{E} := E \setminus (E^0_h \cap A^0)$.
Let $\bar{S}$ be a compact oriented surface 
which is equal to $\pi(\bar{E})$. This has a positive (resp. negative)
boundary component for each positive (resp. negative) puncture of $S$.
The map $\pi : \bar{E} \twoheadrightarrow \bar{S}$
is called a {\it compact convex Lefschetz fibration} associated to $E$.
\end{defn}

This is also defined in \cite[Definition 2.14]{McLean:symhomlef}.
We have that $\bar{E}$ is a manifold with corners.
By \cite[Theorem 2.15]{McLean:symhomlef}, we can assume (maybe after deforming
the Lefschetz fibration without changing the monodromy maps)
that $\bar{E}$ is a Liouville domain after smoothing the corners.
From now on by abuse of notation,  we will write $\bar{E}$
for the Liouville domain obtained by smoothing the corners of $\bar{E}$.
Also if we have a smooth family of compact convex Lefschetz fibrations $\bar{E}_t$,
then the associated Lefschetz fibrations $E_t$ are all symplectomorphic.
Let $\pi : E \twoheadrightarrow \C$ be a Lefschetz fibration and
let $\pi : \bar{E} \twoheadrightarrow \D$ be the associated
compact convex Lefschetz fibration
whose base is the unit disk $\D \subset \C$.
Let $F = \pi^{-1}(1)$.
Let $L$ be an exact Lagrangian sphere in $F$.
This means that $\theta_F|_L = d\kappa$ form some function $\kappa : L \rightarrow \R$
($\theta_F$ is the Liouville form).

We can choose our Liouville form on $\bar{E}$
so that it becomes a Legendrian sphere in $\partial{\bar{E}}$.
The reason is as follows: Because $L$ is an exact Lagrangian in $F$,
we have a small neighborhood of $L$ diffeomorphic to a small
neighborhood of the zero section of $T^*L$ such that
$\theta_F$ is equal to $d\kappa' + \sum_j p_j dq_j$
where $p_j$ are momentum coordinates and $q_j$ are position coordinates.
We can then modify $\theta_F$ so that $\kappa' = 0$ near $L$,
hence $TL$ is in the kernel of $\theta_F$.
The subset $\pi^{-1}(\partial \D) \subset \partial{\bar{E}}$
is a mapping torus of some symplectomorphism $\phi : F \rightarrow F$,
and we can ensure that the contact form is equal
to $Cdt + \theta_F + dG$ as described in section \ref{section:lefschetzfibrations}.
We can also ensure that $G$ is $0$ near $F$.
This then ensures that $L$ is a Legendrian with respect to this contact form.

\begin{theorem} \label{theorem:handleattach}
We can attach a Weinstein $n$ handle $H$ to
$\bar{E}$ along $L$ creating a new Lefschetz fibration
\[ \pi' :\bar{E} \cup H \twoheadrightarrow \D \]
where:
\begin{enumerate}
\item
If $\Delta$ is an arbitrarily small neighborhood of $1 \in \D \subset \C$, then
we can ensure that $\pi'$ coincides with $\pi$ outside $\pi^{-1}(\Delta)$.
\item
$\pi'$ has one extra singularity inside $H$
and its vanishing cycle is Hamiltonian isotopic to $L$ inside $F$.
\end{enumerate}
\end{theorem}

Weinstein handles are described in \cite{Weinstein:contactsurgery} and
\cite{Cieliebak:handleattach}.
We will start by describing the Weinstein $n$-handle $H$:
Let $(p_1,p_2,\cdots,p_n,q_1,\cdots, q_n)$
be standard coordinates for $\C^n$
where the symplectic form $\omega_{\text{std}} = \sum_j dp_j \wedge dq_j$.
We also assume that $z_j = p_j + i q_j$
are the standard complex coordinates form $\C^n$.
Let $p : \C^n \twoheadrightarrow \C$ be defined
by $p(z) = \sum_j z_j^2$.
Let $p_{\R}$ be the real part of $p$,
and $p_{i\R}$ the complex part of $p$.
We write $z = (z_1, \cdots, z_n)$.
We have
$p_{\R}(z) = \sum_j \left( p_i^2 - q_i^2 \right)$
and $p_{\C}(z) = 2 \sum_j \left( p_i q_i \right)$.
Let $x = \sum_j p_j^2$, $y := \sum_j q_j^2$.
We can view $p_\R$ as a function of $x$ and $y$.
This means that $p_\R = x - y$.
From now on let $0 < \delta \ll \epsilon \ll 1$.
We will now construct a function $\psi_{\delta,\epsilon} : \R^2 \rightarrow \R$
which we will also view as a function of $x$ and $y$.
This is constructed as follows:
Let $V_\epsilon$ be a vector field
which is $0$ outside a ball of radius $\epsilon$.
We also assume that $V_\epsilon$ is of the
form $a(x,y) \frac{\partial}{\partial x} + b(x,y) \frac{\partial}{\partial y}$
where $a \leq 0$ and $b \geq 0$ and on the ball of radius $2\delta$ centered at $0$, we set $a = -1,b=1$.
Let $\phi^t_\epsilon : \R^2 \rightarrow \R^2$ be the time $t$
flow of the vector field $V_\epsilon$.
We define $\psi_{\delta,\epsilon} := p_\R \circ \phi^\delta_\epsilon$.
Here is a picture of the level curves of $\psi$ (solid lines) and $p_\R$ (dotted lines):
\begin{figure}[H] 
\centerline{
   \scalebox{1.0}{
 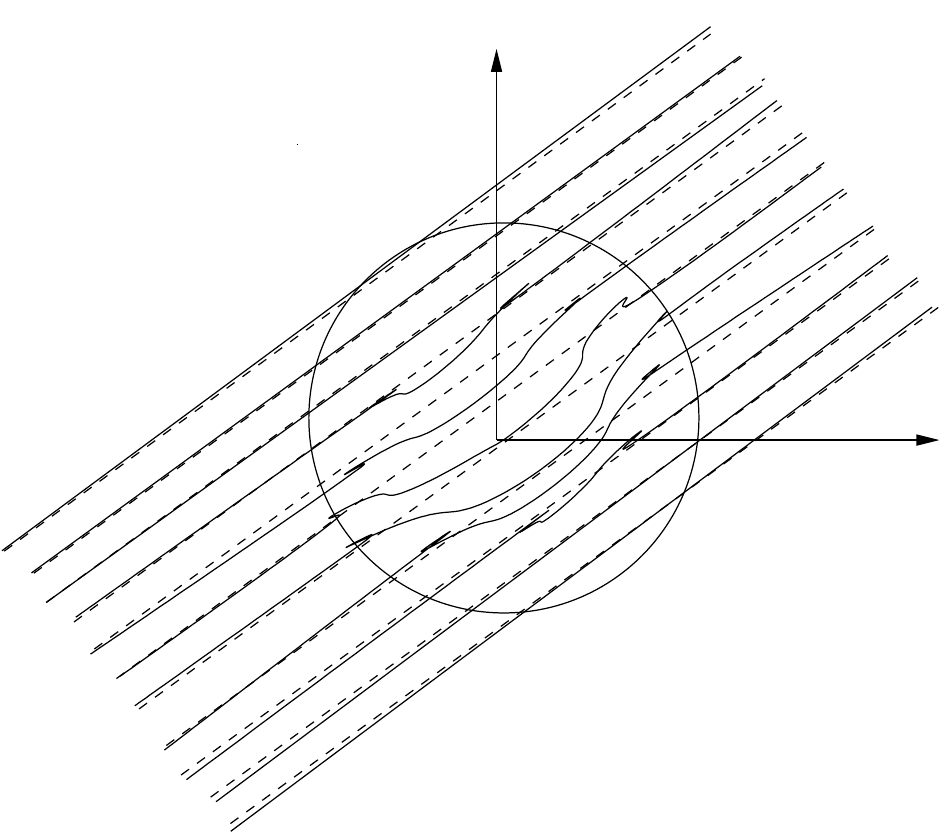
 }
 } 
\end{figure}

We let 
\[H := B_{2\epsilon} \cap \{\psi_{\delta,\epsilon} \leq -\delta/2\} 
\cap \{p_\R \geq -\delta /2 \} \]
where $B_{2\epsilon}$ is a ball of radius $2\epsilon$.
This set is our $n$-handle $H$.
This has a Liouville vector field
\[ \frac{1}{2}\sum_j \left(p_j\frac{\partial}{\partial p_j} -
q_j\frac{\partial}{\partial q_j} \right).\]
We have that $\psi_{\delta,\epsilon}$
is a Weinstein function for this handle
for $0 < \delta \ll \epsilon \ll 1$.
We also define a function
$p' : \C^n \rightarrow \C$ by
$p' = \psi_{\delta,\epsilon} + i p_\C$.
This has a Lefschetz singular point at $0$,
and the map is smooth away from $0$.
Also if $\delta$ is small enough then
the symplectic form $\omega_{\text{std}}$
is non-degenerate on all the vertical tangent
spaces of $p'$ away from $0$.
The attaching region $H_-$ is the subset
$\{p_\R = -\delta / 2\} \cap H \subset H$
and $H_+$ is the subset
$\{\psi_{\delta,\epsilon} = -\delta / 2\} \cap H \subset H$.
Let $S \subset {p^{-1}}(-\delta / 2)$
be the vanishing cycle of the singular point $0$
of $p$ along the path joining $0$ and $-\delta / 2$ along the real axis.
We have that a small neighborhood of $H_-$ in $\left\{ p_\R = -\delta/2 \right\}$ is a codimension $0$
 submanifold of the product
$T^*S^{n-1} \times \R = \left\{ p_\R = -\delta/2 \right\}$ containing
the Legendrian sphere $S \times \left\{ 0 \right\}$.
Here $S$ is the zero section of $T^*S^{n-1}$.
This is true because parallel transport maps for $p$
are well defined, so we can trivialize
the fibration $p_\R|_{\left\{ p_\R = -\delta/2 \right\}}$.
The Liouville vector field
\[ \frac{1}{2}\sum_j \left(p_j\frac{\partial}{\partial p_j} -
q_j\frac{\partial}{\partial q_j} \right)\]
is transverse to $\left\{ p_\R = -\delta/2 \right\}$,
so it is a contact manifold.
The map $p$ restricted to this contact
manifold is smooth and has fibers that are
symplectomorphic to $T^*S^{n-1}$.

We now need a Lemma describing how we can glue
our handle while preserving our fibration structure.
Let $U$ be a contact manifold and $q : U \twoheadrightarrow (-\epsilon',\epsilon') \subset \R$
be a fibration whose fibers are symplectic submanifolds
with symplectic form $d\lambda$ (where $\lambda$ is the contact form). Similarly let
$q' : U' \twoheadrightarrow (-\epsilon',\epsilon')$ be a map with
the same properties as $q$.
We also assume that $q^{-1}(0)$ is exact symplectomorphic
to ${q'}^{-1}(0)$.
Because the fibers are symplectic submanifolds,
there is a natural connection on $q$ and $q'$.
For $q$ this is the line field given by the
kernel of $d\lambda$ and for $q'$
it is the kernel of $d\lambda'$.
This means that we can do parallel transport.
The problem is that it is possible for
a point $p$ to be transported off the edge
of $U$.
We will therefore assume that this cannot
happen for $q$ and $q'$.
Let $t$ parameterize $(-\epsilon,\epsilon)$. We assume that
$\lambda(\widetilde{\frac{\partial}{\partial t}}) > 0$ and
$\lambda'(\widetilde{\frac{\partial}{\partial t}}) >0$
where $\widetilde{\frac{\partial}{\partial t}}$
is the horizontal lift of $\frac{\partial}{\partial t}$
for $q$ or $q'$.

\begin{lemma} \label{lemma:fibreglueing}.
There is a diffeomorphism $C: U \rightarrow U'$
such that $q' \circ C = q$ and such that
$C^* \lambda' = \lambda + dR$ for some function $R$.
Also, there is a smooth family of contact forms $\lambda^s$
such that $d\lambda^s = d\lambda$ for all $s$, $\lambda^0 = \lambda$ and $\lambda^1 = C^*\lambda'$.
\end{lemma}
\proof
%
By using parallel transport maps, we have that
$U$ is diffeomorphic to
$q^{-1}(0) \times (-\epsilon',\epsilon')$
with $\lambda = \lambda|_{q^{-1}(0)} + dR_1$
where $R$ is a function with $\frac{\partial R_1}{\partial t} > 0$.
Similarly we have
$\lambda' = \lambda'|_{{q'}^{-1}(0)} + dR_1'$.
We have a diffeomorphism $\Phi : q^{-1}(0) \rightarrow {q'}^{-1}(0)$
such that $\Phi^*\lambda' = \lambda + dR_2$ for some function $R_2$.
We define $C : U \rightarrow U'$
as $(\Phi, \text{id})$, then
$C^* \lambda' = \lambda - dR_1 + d(R_1' \circ \Phi) + dR_2$.
This proves the first part of the Lemma with $R = -R_1 + R_1' \circ \Phi + R_2$.
We have that $\frac{\partial R_1}{\partial t} > 0$ and
$\frac{\partial (R_1' \circ \Phi + R_2)}{\partial t} > 0$
because $\partial_t R_2 = 0$ and $\partial_t R_1 \circ \Phi >0$.
So we can join these functions via a family of functions $R^s$ with $\frac{\partial R^s}{\partial t} > 0$.
We define $\lambda^s := \lambda|_{q^{-1}(0)} + dR^s$.
These are contact forms because $\lambda|_{q^{-1}(0)}$
is a symplectic form and $\frac{\partial R^s}{\partial t} > 0$,
and $\lambda^0 = \lambda$ and $\lambda^1 = \lambda + dR = C^* \lambda'$.
Also $d\lambda^s = d\lambda$ for all $s$.
This proves the Lemma.
%
%
\qed

Let $f_1 : A_1 \twoheadrightarrow \D$,
$f_2 : A_2 \twoheadrightarrow \D$ be two smooth fibrations
such that $A_1$ and $A_2$ have symplectic forms $\omega_1$
and $\omega_2$ respectively such that the fibers of
$f_1$ and $f_2$ are symplectic submanifolds.
These fibrations have a connection given by
the plane field which is $\omega_i$ orthogonal to the fibers.
We also assume that this connection gives us
well defined parallel transport maps
(i.e. no points get transported off the edge of the manifold $A_i$).
\begin{lemma} \label{lemma:symplecticfibregluing}
Suppose that there exists a line $l \subset \partial \D$
and a fibrewise diffeomorphism $\Phi$ from $A_1$ to $A_2$
such that it is a symplectomorphism from $A_1|_l$ to
$A_2|_l$. We suppose that $\D$ smoothly deformation retracts onto $l$.
Then we can deform $\omega_2$ through symplectic forms $\omega^t_2$
such that $\omega^t_2$ restricted to each fiber is the same as $\omega_2$
and such that there is a symplectic embedding of $A_1$ into some
arbitrarily small neighborhood of $\pi_2^{-1}(l)$ with symplectic form $\omega^1_2$.
This symplectic embedding is equal
to $\Phi$ in the region $\pi_1^{-1}(l)$ and is a fibrewise diffeomorphism.
This deformation also fixes $\omega^t_2$ in the region $\pi_2^{-1}(l)$.
\end{lemma}
\proof of Lemma \ref{lemma:symplecticfibregluing}.
Fix a small neighborhood $(-1,0] \times l \subset \D$
of $l$. The line $l$ is a line with boundary, so we
should lengthen $l$ slightly so that we get a neighborhood.
We identify $l$ with $\left\{ 0 \right\} \times l$ here.
Let $l(0)$ be a point in $l$.
Let $F_i$ be the fiber $\pi_i^{-1}(0)$ for $i=1,2$.
Let $\omega_{F_i}$ ($i=1,2$) be the corresponding symplectic forms
on these fibers.
For a path ${\mathcal P} : [0,1] \twoheadrightarrow \D$, let
$\Phi_{\mathcal P} : \pi_1^{-1}({\mathcal P}(0)) \rightarrow \pi_1^{-1}({\mathcal P}(1))$
be the corresponding parallel transport map for $\pi_1$.
For each point $x = (x_0,x_1) \in (-1,0] \times l$ we have a path $q_x$
joining $l(0)$ with $x$. This path first  travels along $l$
from $l(0)$ to $(0,x_1) \in (-1,0] \times l$ and then
it joins $(0,x_1)$ to $(x_0,x_1)$ along the path $(-1,0] \times \left\{ x_1 \right\}$.
We have the following trivialization:
\[T_1 : (-1,0] \times l \times F_1 \rightarrow \pi_1^{-1}\left( (-1,0] \times l  \right)\]
given by:
$T_1(x_0,x_1,c) = \Phi_{q_{(x_0,x_1)}}(c)$.
The symplectic form on $(-1,0] \times l \times F_1$
is the product $\omega_\D \times \omega_{F_1}$
where $\omega_\D$ is a symplectic form on $(-1,0] \times l$.
We have a similar trivialization $T_2$.
Because $\omega_1$ and $\omega_2$ agree on the region
$\pi_1^{-1}(l) = \pi_2^{-1}(l)$, we have that the symplectic fibrations
are exactly the same via the symplectomorphism $T_2 \circ T_1^{-1}$ (as long
as we choose appropriate neighborhoods of $l$ for $T_1$ and $T_2$ respectively).
Hence we only need to embed the fibration $\pi_1$ into
the fibration given by the image of $T_1$.

Let $r_t : \D \twoheadrightarrow \D$ be the deformation retraction onto $l$.
Here we assume that $r_t$ is a smooth embedding for $t < 1$.
We have that $r_{1-\epsilon}$ maps $\D$ to a small neighborhood of $l$
(so that it fits inside the neighborhood described above.
We have a smooth fibrewise embedding $\Psi$ from
$A_1$ into $\pi_1^{-1}\left( (-1,0] \times l \right)$
given by
$\Phi_{r_t(\pi_1(x))}(x)$.
Here $r_t(\pi_1(x))$ is the path
\[ a : [0,1-\epsilon] \rightarrow \D\]
\[a(t) = r_t(\pi_1(x)).\]
This fibrewise embedding sends fibers to fibers
symplectically. The problem is that it isn't a symplectomorphism.
If we look at $\Psi^*(\omega_1) - \omega_1$, it is a closed
$2$-form which vanishes on the fibers.
Because $\D$ is contractible, this means that
$\omega_1 - \Psi^*(\omega_1) = d\theta$.
Choose a $1$-form $\theta'$ such that $\Psi^*\theta' = \theta$
(we can choose our symplectic embedding $\Psi$ appropriately so that
this works).
Let $h_t : \D \rightarrow \R$ be a family
of functions such that $h_t = t$ near $l$
and is $0$ outside a small neighborhood of $l$.
In order for $\Psi$ to be a symplectic embedding, we
first consider the family of closed $2$-forms
$\omega_1 + d((h_t \circ \pi_1) \theta')$.
All these symplectic forms agree with $\omega_1$
when we restrict to each fiber, but they are not
necessarily symplectic forms on the total space $A_1$.
In order to make them into symplectic forms, we need them to
be non-degenerate. This is done by pulling back
a sufficiently large multiple of a time dependent volume form $v_t$ for $\D$.
This volume form can be chosen so that it is equal to $0$
for $t = 0$ and outside a small neighborhood of $l$.
We can also assume that for $t=1$ it is equal to $0$
on a small neighborhood of $l$ as well.
We have that $\Psi$ is a symplectic embedding if we consider
the symplectic form $\omega_1 + d(h_1 \circ \pi_1) \theta'$.
This proves the lemma.

\qed

\proof of theorem \ref{theorem:handleattach}.
We choose a small neighborhood  $NL \subset F$
of $L$ which is exact symplectomorphic to the
interior of the unit disk bundle $D^*L$
associated to some metric $g$ on $L$.
We also have a similar neighborhood $NS$
of $S \subset p^{-1}(-\delta / 2)$ where
$p : \C^n \twoheadrightarrow \C$ is the map
described earlier.
Choose  $NL$ and $NS$ so that they are exact symplectomorphic.
The exponential map gives us a smooth embedding
$e : (-\frac{1}{2},\frac{1}{2}) \rightarrow \partial \D$.
Let $P_t$ be the parallel transport map for $\pi$
starting at $e(0)$ and traveling along the path
$e$ to $e(t)$.
We define $U := \cup_{|t|<\frac{1}{2}} P_t(NL)$.
We let $e' : (-\frac{1}{2},\frac{1}{2}) \rightarrow \C$
be the path defined by $e'(t) = -\delta / 2 + it$.
Let $P'_t$ be the parallel transport map for $p$
starting at $e'(0)$ and traveling along the path
$e'$ to $e'(t)$.
We define $V$ to be $\cup_{|t|<\frac{1}{2}} P'_t(NS) \subset \C^n$.
By using Lemma \ref{lemma:fibreglueing} we have a
diffeomorphism $C : U \rightarrow  V$ such that
$C^* p|_V = \pi|_U$ and such that $C^*d\theta_H|_V = d\theta_E|_U$
where $\theta_H$ is the Liouville form on the Handle and
$\theta_E$ is the Liouville form on the Lefschetz fibration $\bar{E}$. Also by the second part of this
Lemma, we can deform the contact form on the boundary of $\bar{E}$
through contact forms so that $C$ becomes a contactomorphism.
We can ensure that this deformation of contact forms extends
to a Liouville deformation of $\bar{E}$.
The handle $H$ depends on parameters $\delta$ and $\epsilon$,
so for $0 < \delta \ll \epsilon \ll 1$, we have that
$H_-$ is a subset of $V$.
Hence we can use $C^{-1}$ to glue $H_-$
to $U$.
We will write $\bar{E} \cup H$ for the Lefschetz fibration with the handle glued.
We now wish to extend the Lefschetz fibration structure
over the handle $H$.
The handle $H$ has two maps
$p$ and $p'$ from $H$ to $\C$.
We will now enlarge $H$ to $\widetilde{H}$
inside $\C^n$ as follows:
We use parallel transport to construct
a subset $[-\delta,-\delta/2] \times H_-$
where the map $p_\R$ is the projection
to $[-\delta,-\delta/2]$
and the map $p_\C$ is the composition of 
the projection to $H_-$ with
$p_\C|_{H_-}$.

We can use the same parallel transport trick
to construct a neighborhood of $U$
in $E$ diffeomorphic to $[-\delta,-\delta/2] \times U$
as follows:
We have a small neighborhood of the image
of $e$ inside the disk $\D$ is biholomorphic
to the product $[-\delta,-\delta/2] \times (-\frac{1}{2},\frac{1}{2})$.
We can ensure that it is disjoint from the singular values of $\pi$.
Then we use parallel transport as before
to construct a subset of $E$ diffeomorphic
to $[-\delta,-\delta/2] \times U$
which projects via $\pi$ to $[-\delta,-\delta/2] \times (-\frac{1}{2},\frac{1}{2})$.
We have a smooth fibrewise embedding of $[-\delta,-\delta/2] \times H_-$
into $[-\delta,-\delta/2] \times U$.
Each inclusion $\left\{ x \right\} \times H_- \hookrightarrow \left\{ x \right\} \times U$
is a symplectic embedding.
By Lemma \ref{lemma:symplecticfibregluing}, we can Liouville deform $\bar{E}$ through compact convex
Lefschetz fibrations so that it becomes
a symplectic embedding.
The map $\pi$ coincides with $p$
in this region.
Let $N \subset H_-$ be a neighborhood of the boundary of $H_-$ inside $H_-$.
Let $H^U$ be the region
\[\left(  \left(  [-\delta,-\delta+\eta] \right) \times H_- \right) 
\cup \left(  [-\delta,-\delta/2] \times N \right) \subset H.\]
For $\eta$ small enough, the map $p'$ coincides with $p$
in the region $F$.
This means we can change $\pi$ to a new map $\pi'$
in the following way:
We define $\pi'$ to be equal to $\pi$ outside
the region $[-\delta,-\delta/2] \times H_-$ and
set it equal to $p'$ inside this region.
This map is smooth because $p'$ coincides
with $\pi$ on the region $F$.
The map $\pi'$ also has one extra singularity
whose vanishing cycle is Lagrangian isotopic
to the original Lagrangian sphere $L$.
We also have that $\bar{E} \cup H$ is equal to $\bar{E}$
with an $n$-handle attached
\qed

\subsection{Showing handle cancellation} \label{section:appendixmainargument}

Let $F$ be a Liouville domain of dimension $2n -2$.
Let $\phi : \widehat{F} \rightarrow \widehat{F}$ be a compactly supported symplectomorphism
where $\widehat{F}$ is the completion of the Liouville domain $F$.
We have that $\left\{ r_F \leq K \right\}$ is also a Liouville domain
for any $K$ whose completion is also symplectomorphic to $\widehat{F}$.
The support of $\phi$ is contained in this set for large enough $K$.
Because of this we may as well redefine $F$ so that the
support of $\phi$ is contained inside $F$.
Let $F'$ be obtained from $F$ by attaching
an $n-1$-handle.
The $n-1$-handle $H$ has a natural Weinstein function
given by $\psi$ described in the previous section.
This has exactly one critical point.
This has an unstable manifold.
Let $L$ be a Lagrangian in $F'$ which intersects the unstable
manifold in exactly one point.
Let $\tau_L$ be a symplectomorphism which is a Dehn twist
around $L$.
Because $\phi$ is the identity near the boundary of $F$,
we can extend $\phi$ to a symplectomorphism  $\phi'$ from $F'$
to $F'$ by making $\phi$ equal to the identity on the handle.
A Stabilization of $\phi$ is defined to be
$\tau_L \circ \phi$.
We can extend this map to the completion $\widehat{F'}$
of $F'$ by making it equal to the identity map outside $F'$.
By abuse of notation, we will use the same names
$\tau_L$ and $\phi$ for the corresponding maps defined on the completion.

\begin{theorem} \label{thm:stabilizationtheorem}
Let $\pi : \widehat{E} \twoheadrightarrow \C$
be a Lefschetz fibration. Let $\phi : \widehat{F} \rightarrow \widehat{F}$ be the monodromy
map.  Then there exists another Lefschetz fibration
$\pi'' : \widehat{E''} \twoheadrightarrow \C$ whose monodromy map
is $\tau_L \circ \phi : \widehat{F'} \rightarrow \widehat{F'}$
such that $\widehat{E''}$ is symplectomorphic to $\widehat{E}$.
\end{theorem}

\proof of Theorem \ref{thm:stabilizationtheorem}.
The Lefschetz fibration $\widehat{E}$ is the completion
of some compact convex Lefschetz fibration $E$
see \cite[Definition 2.16]{McLean:symhomlef}).
The fiber of this Lefschetz fibration is some Weinstein domain $F$.
A small neighborhood of $\partial F$ is symplectomorphic to
$(1-\eta, 1] \times \partial F$ with Liouville form $r_F \alpha_F$
where $\alpha_F$ is the contact form on the boundary of $F$ and $r_F$
is the coordinate parameterizing $(1-\eta, 1]$.
Because $E$ is a Liouville domain (with corners),
we have a Liouville vector field $X$ on $E$ which is transverse to the boundary of $E$.
Its boundary is a union of two manifolds $E^h$ and $E^v$ meeting in one corner $\{\partial E^h \cap \partial E^v\}$.
The manifold $E^v$ is the vertical boundary equal to $\pi^{-1}(\partial \D)$ and
$E^h$ is the horizontal boundary.
A small neighborhood of $E^h$ is symplectomorphic to
$\D \times \left( (1-\eta,1] \times \partial F \right)$.
We can ensure that the Liouville vector field $X$ is equal to
$r \frac{\partial}{\partial r} + r_F \alpha$ in this region where
$r$ is the radial coordinate for $\D$.

A small neighborhood of $\partial{E}$ is of the form
$\left\{ r_F \geq 1-\eta \right\} \cup \left\{ r \geq 1- \eta \right\}$.
Let $g: (1-\eta,1] \rightarrow \R$ be a function which is equal to
$0$ on a small neighborhood of $1-\eta$ and such that $g(x) = x$
for $x$ near $1$.
We can define the function $g(r_F) + r^2 : E \rightarrow \R$.
Because $g(x)$ is equal to zero near $1-\eta$,
this function is well defined on all of $E$ as we extend it by $r^2$
when $r_F$ is ill defined.
Here we write $r$ instead of $\pi^* r$ by abuse of notation.
This is also a Weinstein function with respect to $X$
on a neighborhood of $\partial E$ (but not on all of $E$).

We then attach the $n-1$ handle to $F$ as above to create
a new Liouville domain $F'$. Let $U$ be a small neighborhood
of $F' \setminus F$. We have a Weinstein function $w : U \rightarrow \R$
with exactly one singular point of index $n-1$ and such that
$w = r_F$ near $\partial F \subset F'$.
Let $\bar{H}$ be this handle. This is equal to the closure
of $F' \setminus F$ inside $F'$.
Let $Y$ be the Liouville vector field described above for this handle.
This is has one critical point of index $n-1$ and $w$
is a Weinstein function for this handle.

Let $E'$ be a new compact convex Lefschetz fibration
obtained by gluing $\D \times \bar{H}$ to the region
$E^h = \D \times \partial F$ where we identify
the $\D$ factor  with itself and we attach $\bar{H}$
to $\partial F$.
The map $\pi$ extends to a map $\pi' : E' \twoheadrightarrow \D$
by setting $\pi'$ equal to the projection map
\[\D \times \bar{H} \twoheadrightarrow \D\] outside $E$.
We can extend the function $g(r_F) + r^2$
to $E'$ over the region $\D \times \bar{H}$
by the function $r^2 + w$.
Let $w_{E'} : E' \twoheadrightarrow \R$ be this new function.
Also, the Liouville vector field
$X$ extends to a Liouville vector field $X'$ defined on $E'$ by setting
$X'$ to be equal to
\[r \frac{\partial}{\partial r} + Y\]
on the region $\D \times \bar{H}$
and $X' = X$ elsewhere
The vector field $X'$ has exactly
one critical point of index $n-1$ in the region $E' \setminus E$.
The unstable manifold of this critical point is of the
form $\D \times U_Y$ inside $E' \setminus E$ where $U_Y$ is the unstable manifold of $Y$.

We have a Lagrangian $L \subset F'$.
Without loss of generality, we can assume
that $0$ is a regular value of $\pi'$.
Identify $F'$ with ${\pi'}^{-1}(0)$.
We now attach our $n$ handle as in
Theorem \ref{theorem:handleattach} along $L$.
This means we need to modify the Liouville domain $E'$
on $\pi^{-1}(U)$ where $U$ is an arbitrarily small neighborhood of $1$.
We can assume that $U$ is disjoint from all the singular values and also $0$.
Let $\pi'' : E'' := E' \cup H \twoheadrightarrow \D$ be this new Liouville domain, and let $X''$
be the new Liouville vector field on $E''$.
We can extend the Weinstein function $w_{E'}$ over this handle
so it becomes a Weinstein function $w_{E''}$ for $X''$.
The unstable manifold $\D \times U_Y$ described above
is also an unstable manifold for the same critical point of $X''$ away
from $\pi^{-1}(U)$. Let $A \subset E''$ be the unstable manifold of $X''$
at this critical point.
The Liouville vector field $X''$ has one extra critical point
in the handle of index $n$.
The map $\pi''$ has one extra critical value $x \in \D$ near $1$.
Let $l$ be a path joining $x$ with $0$ avoiding all the critical
values of $\pi''$.
Let $V \subset E''$ be the thimble associated to this path.
I.e. it is the set of points in $E''$ that parallel
transport along this path into the singularity of $\pi''$ inside the handle.
This is a Lagrangian submanifold diffeomorphic to a
ball (see \cite[Section 1]{Seidel:longexactsequence}).
Let $V_1 \subset E''$ be the stable manifold of the singularity
of $X''$ in the handle $H$.
Near the singularity $x$ we have that $V_1$ is the
same as $V$.
This is because we have that $\pi''$ is the same (after a translation)
as $p : \C^n \rightarrow \C$. By symmetry (the antiholomorphic
involution fixing $\R^n \subset \C^n$)
we have that $V$ and $V_1$ are open subsets
of $\R^n \subset \C^n$.
We also have for some $\epsilon > 0$ small enough,
that $\tau := w_{E''}^{-1}(w_{E''}(x) - \epsilon)$
is a contact submanifold which is isotopic through contact
submanifolds to some smoothing of the boundary
of $E' \subset E''$.
Also because the Lagrangian sphere $(\pi'')^{-1}(l(0)) \cap V$
is Hamiltonian isotopic inside $(\pi'')^{-1}(l(0))$ to a Lagrangian sphere intersecting
$A$ once, we have that the corresponding Legendrian sphere
$\partial E' \cap V$ is Legendrian isotopic to a Legendrian sphere
intersecting $A$ once.
Hence $\tau$ is Legendrian isotopic to a Legendrian sphere intersecting
$A$ once.

By work of Eliashberg \cite[Lemma 3.6 b]{Eliashberg:symplecticgeometryofplushfns}
we can replace the Weinstein function $w_{E''}$ and the Liouville form
with another Weinstein function $w'_{E''}$ such that
$w'_{E''}|_E = w_{E''}|_E$ and such that $w'_{E''}$
has no critical points outside $E$.
This ensures that the completion $\widehat{E}$
is symplectomorphic to $\widehat{E''}$.
\qed

%

\section{Appendix B: A minimum principle}

We have a Lefschetz fibration $\pi : E \rightarrow \C$.
Let $[1,\infty) \times M_\phi$ be its cylindrical end where $r_S$ parameterizes the interval.
The $1$-form $\theta_E$ is equal to $r_s d\vartheta + \alpha_\phi$
in this region where $\alpha_\phi$ is the contact form on $M_\phi$
and $\vartheta$ is the pullback of the angle coordinate on $S^1$ to $M_\phi$.
Let $H_{s,t}$ be a family of Hamiltonians parameterized by $\R \times S^1$ where
$H_{s,t} = \kappa r_S + g(r_F)$ where $\kappa$ is a constant near the level set $r_S = \Delta$
and $r_F$ is the fiber cylindrical coordinate.
Let $J_t$ be an $S^1$ family of almost complex structures where
$\pi$ is $(J_t,j)$ holomorphic near $r_S = \Delta$.
Here $j$ is the standard complex structure on the cylinder $[1,\infty) \times S^1$.
\begin{lemma} \label{lemma:lefschetzminimumprinciple}
Any Floer trajectory which does not intersect a regular fiber $\pi^{-1}(q)$
connecting orbits in the region $r_S > \Delta$
must be contained in the region $r_S > \Delta$.
\end{lemma}
\proof of Lemma \ref{lemma:lefschetzminimumprinciple}.
The proof of this is very similar to the proof of \cite[Lemma 7.2]{SeidelAbouzaid:viterbo}.
A Floer trajectory connecting orbits is a cylinder:
\[u : \R \times S^1 \rightarrow \R \text{, }
u_s + J_t u_t = J_t X_{H_{s,t}}.\]
After perturbing $\Delta$ slightly we have that $u$ is transverse to $\{r = \Delta\}$
so $u^{-1}(\{r_S \leq \Delta\})$ is a compact codimension $0$ submanifold $\overline{S}$
of the cylinder  $\R \times S^1$.
There is a $1$ form $\beta$ on
$E \setminus \pi^{-1}(q)$ such that
$\beta|_{M_\phi} = d\vartheta$ where $\vartheta$ is the angle
coordinate for $S^1$.
The Hamiltonian vector field associated to $H_{s,t}$ near $r = \Delta$ is equal to
$-\kappa \widetilde{\frac{\partial}{\partial \theta}} + X$
where $\widetilde{\frac{\partial}{\partial \theta}}$ is the horizontal lift of
$\frac{\partial}{\partial \theta}$ and $X$ is the Hamiltonian flow of $g(r_F)$.
Here $X$ is tangent to the fibers of $\pi$.
We have:
\[ \int_{\partial \overline{S}} u^* d\vartheta = \int_{\partial \overline{S}} u^*\beta
= \int_{\overline{S}} u^* d\beta = 0.\]
Hence
\[ 0 = \int_{\partial \overline{S}} \kappa dt =
\int_{\partial \overline{S}} u^* d\vartheta - d\vartheta \left(-\kappa \widetilde{\frac{\partial}{\partial \theta}} + X\right) dt
=\int_{\partial \overline{S}} u^* d\vartheta - d\vartheta X_{H_{s,t}} dt\]
\[ = \int_{\partial \overline{S}} d\vartheta(u_s) ds + d\vartheta(u_t - X_{H_{s,t}}) d
 = \int_{\partial \overline{S}} d\vartheta(J_t(X_{H_{s,t}})-J_t u_t) ds + d\vartheta(J_t u_s) dt\]
Now $d\vartheta \circ J_t = dr_S$ and $dr_S(X_{H_{s,t}}) = 0$ along $\{r_S = \Delta\}$
hence our integral becomes:
\[ \int_{\partial \overline{S}} dr_S( u_s) dt + dr_S(X_{H_{s,t}} - u_t) ds 
= \int_{\partial \overline{S}} dr_S(u_s) dt - dr_S(u_t) ds\]
\[= \int_{\partial \overline{S}} dr_S \circ du \circ j.\]
Let $\xi$ be a vector on $\partial \overline{S}$ which is positively oriented then
$j(\xi)$ points inwards along $\partial \overline{S}$ and because $u$ is transverse to $\{r = \Delta\}$
we get that $dr_S \circ du \circ j(\xi) < 0$ which implies that our integral satisfies:
\[\int_{\partial \overline{S}} dr_S \circ du \circ j < 0\] which is  a contradiction.
\qed

\bibliography{references}

\newcommand{\etalchar}[1]{$^{#1}$}
\begin{thebibliography}{CFHW96}

\bibitem[{Akb}10]{Akbulut:lefschetzfibrationsoncompactstein}
{Akbulut, S}.
\newblock {Lefschetz Fibrations on Compact Stein Manifolds}.
\newblock pages 1--19, 2010, arXiv:1003.2200.

\bibitem[AM09]{McLeanAlbers:leafwise}
P.~Albers and M.~McLean.
\newblock Non-displaceable contact embeddings and infinitely many leaf-wise
  intersections.
\newblock pages 1--11, 2009, arXiv:0904.3564.

\bibitem[AS07]{SeidelAbouzaid:viterbo}
M.~Abouzaid and P.~Seidel.
\newblock An open string analogue of {V}iterbo functoriality.
\newblock pages 1--74, 2007, arXiv:0712.3177.

\bibitem[BEE]{BEE:legendriansurgery}
F.~Bourgeois, T.~Ekholm, and Y.~Eliashberg.
\newblock Symplectic homology as {H}ochschild homology.
\newblock arXiv:SG/0609037.

\bibitem[BEH{\etalchar{+}}03]{BEHWZ:compactnessfieldtheory}
F.~Bourgeois, Y.~Eliashberg, H.~Hofer, K.~Wysocki, and E.~Zehnder.
\newblock Compactness results in symplectic field theory.
\newblock {\em Geom.Topol.}, 7:799--888, 2003, arXiv:SG/0308183.

\bibitem[BO09]{BourgeoisOancea:exactsequence}
F.~Bourgeois and A.~Oancea.
\newblock An exact sequence for contact- and symplectic homology.
\newblock {\em Invent. Math.}, 175(3):611--680, 2009.

\bibitem[CC09]{CottonClay:symplecticfloerarea}
Andrew Cotton-Clay.
\newblock Symplectic {F}loer homology of area-preserving surface
  diffeomorphisms.
\newblock {\em Geom. Topol.}, 13(5):2619--2674, 2009.

\bibitem[CFHW96]{CieliebakFloerHoferWysocki:SymhomIIApplications}
K.~Cieliebak, A.~Floer, H.~Hofer, and K.~Wysocki.
\newblock Applications of symplectic homology {II}:stability of the action
  spectrum.
\newblock {\em Math. Z}, 223:27--45, 1996.

\bibitem[Cie02]{Cieliebak:handleattach}
K.~Cieliebak.
\newblock Handle attaching in symplectic homology and the chord conjecture.
\newblock {\em J. Eur. Math. Soc. (JEMS)}, 4:115--142, 2002.

\bibitem[DS94]{DostoglouSalamon:instantonscurves}
S.~Dostoglou and D.~Salamon.
\newblock Self-dual instantons and holomorphic curves.
\newblock {\em Ann. of Math. (2)}, 139:581--640, 1994.

\bibitem[Eli97]{Eliashberg:symplecticgeometryofplushfns}
Y.~Eliashberg.
\newblock Symplectic geometry of plurisubharmonic functions, notes by {M}.
  {A}breu, in: Gauge theory and symplectic geometry ({M}ontreal 1995).
\newblock {\em NATO Adv. Sci. Inst. Ser. C Math. Phys. Sci.}, 488:49--67, 1997.

\bibitem[FH94]{FloerHofer:SymhomI}
A.~Floer and H.~Hofer.
\newblock Symplectic homology {I}.
\newblock {\em Math. Z}, 215:37--88, 1994.

\bibitem[Gir02]{Giroux:openbooks}
Emmanuel Giroux.
\newblock G\'eom\'etrie de contact: de la dimension trois vers les dimensions
  sup\'erieures.
\newblock In {\em Proceedings of the {I}nternational {C}ongress of
  {M}athematicians, {V}ol. {II} ({B}eijing, 2002)}, pages 405--414, Beijing,
  2002. Higher Ed. Press.

\bibitem[Hat02]{Hatcher:algebraictopology}
A.~Hatcher.
\newblock {\em Algebraic topology}.
\newblock Cambridge University Press, Cambridge, 2002.

\bibitem[HS95]{HoferSalamon:FloerNovikov}
H.~Hofer and D.~A. Salamon.
\newblock Floer homology and {N}ovikov rings.
\newblock In {\em The {F}loer memorial volume}, volume 133 of {\em Progr.
  Math.}, pages 483--524. Birkh\"auser, Basel, 1995.

\bibitem[McL]{McLean:computability}
M.~McLean.
\newblock {Computability and the growth rate of symplectic homology}.
\newblock {\em In preparation}.

\bibitem[McL09]{McLean:symhomlef}
M.~McLean.
\newblock Lefschetz fibrations and symplectic homology.
\newblock {\em Geom. Topol.}, 13(4):1877--1944, 2009.

\bibitem[Oan06]{Oancea:kunneth}
A.~Oancea.
\newblock The {K}{\"u}nneth formula in {F}loer homology for manifolds with
  restricted contact type boundary.
\newblock {\em Math. Ann.}, 334:65--89, 2006, arXiv:SG/0403376.

\bibitem[RS93]{RobbinSalamon:maslov}
J.~Robbin and D.~Salamon.
\newblock The {M}aslov index for paths.
\newblock {\em Topology}, 32:827--844, 1993.

\bibitem[Sei03]{Seidel:longexactsequence}
P.~Seidel.
\newblock A long exact sequence for symplectic {F}loer cohomology.
\newblock {\em Topology}, 42:1003--1063, 2003, arXiv:SG/0105186.

\bibitem[Sei08]{Seidel:biasedview}
P.~Seidel.
\newblock A biased view of symplectic cohomology.
\newblock {\em Current Developments in Mathematics}, 2006:211--253, 2008.

\bibitem[Vit99]{Viterbo:functorsandcomputations}
C.~Viterbo.
\newblock Functors and computations in {F}loer homology with applications, part
  {I}.
\newblock {\em Geom. Funct. Anal.}, 9:985--1033, 1999.

\bibitem[Wei91]{Weinstein:contactsurgery}
A.~Weinstein.
\newblock Contact surgery and symplectic handlebodies.
\newblock {\em Hokkaido Math. J.}, 20(2):241--251, 1991.

\end{thebibliography}

\end{document}